\documentclass{amsart}
\usepackage{amsmath}
\usepackage{amssymb}
\usepackage{extarrows}
\usepackage{mathabx}
\usepackage[all]{xy}%commutative graph
\usepackage{amsbsy}
\usepackage{graphicx}
\usepackage{appendix}
\usepackage{float}
\usepackage{subfigure}
\usepackage[numbers,sort&compress]{natbib}
\usepackage{graphicx}%pictures
\usepackage{amsthm}
\usepackage{geometry}
\usepackage{fancyhdr}%Heads
\usepackage{color} %colors
\usepackage{overpic}%pictures with text
\usepackage{mathabx}
\usepackage{tikz-cd}
\usepackage{enumerate}
\usepackage{mathrsfs}
\usepackage{colonequals}
\usepackage{microtype}
\usepackage{hyperref}
\usepackage{diagbox}

\newtheorem{thm}{Theorem}[section]
\newtheorem{cor}[thm]{Corollary}
\newtheorem{prop}[thm]{Proposition}
\newtheorem{lem}[thm]{Lemma}

\theoremstyle{definition}
\newtheorem{defn}[thm]{Definition}
\newtheorem{cons}[thm]{Construction}
\newtheorem{exmp}[thm]{Example}

\newtheorem*{fact}{Fact}

\theoremstyle{remark}
\newtheorem{rem}[thm]{Remark}

\numberwithin{equation}{section}

\newcommand{\beq}{\begin{equation*}\begin{aligned}}
\newcommand{\eeq}{\end{aligned}\end{equation*}}
\newcommand{\bpf}{\begin{proof}}
\newcommand{\epf}{\end{proof}}
\newcommand{\bthm}{\begin{thm}}
\newcommand{\ethm}{\end{thm}}
\newcommand{\bprop}{\begin{prop}}
\newcommand{\eprop}{\end{prop}}
\newcommand{\bcor}{\begin{cor}}
\newcommand{\ecor}{\end{cor}}
\newcommand{\blem}{\begin{lem}}
\newcommand{\elem}{\end{lem}}
\newcommand{\bdefn}{\begin{defn}}
\newcommand{\edefn}{\end{defn}}
\newcommand{\bcons}{\begin{cons}}
\newcommand{\econs}{\end{cons}}
\newcommand{\bexmp}{\begin{exmp}}
\newcommand{\eexmp}{\end{exmp}}
\newcommand{\brem}{\begin{rem}}
\newcommand{\erem}{\end{rem}}
\newcommand{\bfa}{\begin{fact}}
\newcommand{\efa}{\end{fact}}
\newcommand{\benu}{\begin{enumerate}[(1)]}
\newcommand{\eenu}{\end{enumerate}}

\newcommand{\bdia}{\begin{displaymath}\xymatrix}
\newcommand{\edia}{\end{displaymath}}

\newcommand{\al}{\alpha}
\newcommand{\be}{\beta}

\newcommand{\Ga}{\Gamma}

\newcommand{\dash}{\textrm{-}}%en dash

\newcommand{\intg}{\mathbb{Z}}

\newcommand{\ra}{\rightarrow}

\newcommand{\xra}{\xrightarrow}

\DeclareMathOperator{\im}{Im}

\begin{document}

\title{Knot cobordism, torsion order and framed instanton homology}

%    Remove any unused author tags.

%    author one information

\author{Sudipta Ghosh}
\address{Department of Mathematics, University of Notre Dame}
\curraddr{}
\email{sghosh7@nd.edu}
\thanks{}

\author{Zhenkun Li}
\address{Department of Mathematics and Statistics, University of South Florida}
\curraddr{}
\email{zhenkun@usf.edu}
\thanks{}

%author two information
%\author{Fan Ye}
%\address{Department of Mathematics, Harvard University}
%\curraddr{}
%\email{fanye@math.harvard.edu}
%\thanks{}

%author two information
%\author{Fan Ye}
%\address{Department of Pure Mathematics and Mathematical Statistics, University of Cambridge}
%\curraddr{}
%\email{fy260@cam.ac.uk}
%\thanks{}

\date{}
\dedicatory{}
\maketitle
\begin{abstract}
	We construct cobordism maps for the \textit{minus} version of instanton knot homology associated to a \textit{specially decorated} knot cobordisms of arbitrary genus between two null-homologous knots in closed oriented $3$-manifolds. As an application of our construction, we recover an inequality between the torsion order of knots in instanton theory, which was originally established in Heegaard Floer theory by work of Juh\'asz, Miller, and Zemke. We further use this inequality to compute the framed instanton Floer homology of any non-zero Dehn surgeries along an alternating knot of bridge index at most $3$.
\end{abstract}

\section{Introductions}
In Heegaard Floer homology, link cobordism maps associated to decorated link cobordisms were constructed by Juh\'asz \cite{Juhasz2012} for hat version and by Zemke \cite{zemke2019link} for general version. The construction of link cobordism maps has led to many interesting topological applications, such as the stabilization distance between minimal-genus surfaces in $B^4$ bounding knots in $S^3=\partial B^3$ by Juh\'asz and Zemke \cite{juhasz2018stabilization} and a bound on the bridge index of a knot from the torsion order of the knot by Juh\'asz, Miller, and Zemke \cite{juhasz2020bridge}. It is natural to ask whether similar constructions can be made in other branches of Floer theory. In instanton theory, the instanton knot homology $KHI$ was constructed by Kronheimer and Mrowka \cite{kronheimer2010knots}, which corresponds to the hat version of knot Floer homology. The minus version of instanton knot homology $KHI^-$ was constructed for knots by the authors of the current paper in \cite{li2019direct} and \cite{GL2023decompose}, and was recently revisited by the first author and Zemke \cite{GZ2023}. As for link cobordism maps, the construction for $KHI$ follows from the work of the second author in \cite{li2018gluing}, while the construction for $KHI^-$ was previously unknown.

 Suppose we have a link cobordism $(W,\Sigma):(Y_0,K_0)\to (Y_1,K_1)$. In Zemke's construction, $(W,\Sigma)$ was equipped with some extra data called decorations and is decomposed into a few fundamental types: the birth, death, and saddle. Then the link cobordism map associated to the decorated cobordism is the composition of ones associated to the fundamental types. In this paper we develop a different construction. Instead of decomposing, we use the surface $\Sigma$ and a special decoration $\mathcal{D}$ on $\Sigma$ to construct a map
 $$\hat{\rho}: KHI^-(Y_0,K_0)\to KHI^-(Y_0\#Z_{2g},K_1\# B_g),$$
 where $B_g\subset Z_{2g}$ is the borromean knot as in Section \ref{subsec: borromean knots}, with $g=g(\Sigma)$. Then the $4$-dimensional cobordism $W$ induces a map 
 $$F^-_{W-N(\Sigma)}:KHI^-(Y_0\#Z_{2g},K_0\# B_g)\to KHI^-(Y_1,K_1).$$
 Thus the composition gives rise to the desired cobordism map
 $$F_{W,\Sigma_g,\mathcal{D}}=F^-_{W-N(\Sigma)}\circ \hat{\rho}: KHI^-(Y_0,K_0)\to KHI^-(Y_1,K_1).$$

\begin{thm}\label{thm: cobordism map, intro}
    If $(W, \Sigma_g,\mathcal{D}): (Y_0, K_0) \to (Y_1, K_1)$ is a specially decorated knot cobordism, then the map $F_{W,\Sigma_g,\mathcal{D}}$ constructed as above is $\mathbb {C}[U]$-equivariant and is functorial under the composition of such cobordisms.
\end{thm}

Using the same strategy we also extend our construction from knots to 2-component links $L=K \cup U$ in a disconnected manifold where $K \subset Y$ is an arbitrary knot and $U \subset Y'$ is the unknot. We prove that the knot cobordism maps and these special classes of link cobordism maps are functorial with respect to compositions.

\begin{rem}
	In order to define the map $\hat{\rho}$ and hence $F_{W,\Sigma_g,\mathcal{D}}$, one needs to choose an element in $KHI^-(Z_{2g},B_{2g})$. This homology group has been computed explicitly, and in this paper, in order to establish our applications, we choose a special element in $KHI^-(Z_{2g},B_{2g})$ and omit it from the data. However, there are many other elements that could potentially lead to different link cobordisms. So one could potentially choose different elements to define different link cobordism maps.
\end{rem}

\begin{rem}\label{rem: scalar ambiguity}
	Due to the naturality issue ({\it c.f.} \cite{baldwin2015naturality}), all instanton Floer homology groups and maps between them in this paper are only well defined up to the multiplication of a non-zero element in $\mathbb{C}$.
\end{rem}

Instanton Floer homology was constructed based on a set of partial differential equations and is known to be very hard to compute. In recent years, there have been several works towards computing the framed instanton Floer homology of various families of 3-manifolds, especially those coming from Dehn surgeries on knots. See, for example \cite{lidman2020framed,baldwin2021concordance,baldwin2022concordance2,alfieri2020framed,LY2021large,LY2022integral2}. As an application of our construction for the cobordism maps for $KHI^-$, we include a large family of 3-manifolds into the known list.
\bthm\label{thm: main}
Suppose $K\subset S^3$ is an alternating knot with a of bridge index at most 3 and $r\in\mathbb{Q}$ is non-zero. Then we have
$$\dim_{\mathbb{C}}I^{\sharp}(S^3_r(K))={\rm rk}_{\mathbb{Z}_2}\widehat{HF}(S^3_r(K)).$$
\ethm
\brem
Although we state the theorem as an identity between the dimensions in instanton Floer theory and Heegaard Floer theory, given an (alternating) diagram of any such knot, there is an algorithm to compute the precise dimension from the diagram. Also, it should be noted that 2-bridge knots are automatically alternating.
\erem

%Instanton Floer homology is closely related to the representation varieties of fundamental groups, and has many important applications to $3$-manifold topology. However, it is impossible to compute the instanton Floer homology from the definition, and the relation between instanton Floer homology and other Floer theories remains elusive. Kronheimer and Mrowka \cite{kronheimer2010knots} conjectured an isomorphism between the framed instanton Floer homology and the hat version of Heegaard Floer homology, but currently not much is known beyond computational examples.

The proof of Theorem \ref{thm: main} relies on a recent progress in developing the surgery formula in instanton Floer homology by the second author and Ye. In particular, in \cite{LY2022integral2}, we have the following.
\begin{prop}\label{prop: dehn surgery on torsion order one}
	Suppose $K\subset S^3$ is an alternating knot of torsion-order-one, and $r\in\mathbb{Q}$ is non-zero. Then we have
$$\dim_{\mathbb{C}}I^{\sharp}(S^3_r(K))={\rm rk}_{\mathbb{Z}_2}\widehat{HF}(S^3_r(K)).$$
\end{prop}

Regarding the torsion order, the second author constructed in \cite{li2019direct} a minus version of instanton knot homology, which we denote by $KHI^-(S^3,K)$ for a knot $K\subset S^3$. The minus version $KHI^-(S^3,K)$ admits an action of $U$ and the torsion order of $K$ is defined as the maximal $U$-order of any $U$-torsion element in $KHI^-(S^3,K)$. We denote the torsion order of $K$ by ${\rm ord}_{U}(K)$

Utilizing Proposition \ref{prop: dehn surgery on torsion order one} and the fact that the torsion order of alternating knots are odd ({\it c.f.} Lemma \ref{lem: odd torsion order}), in order to prove Theorem \ref{thm: main}, it remains to prove the following.
\bthm\label{thm: torsion order and bridge index}
Suppose $K\subset S^3$ is a knot. Then
$${\rm ord}_{U}(K)\leq {\rm br}(K)-1.$$
Here ${\rm br}(\cdot)$ denotes the bridge index of $K$.
\ethm

A parallel result in Heegaard Floer theory has been established by Juh\'asz, Miller, and Zemke \cite{juhasz2020bridge}. Previously three ingredients were missing to establish Theorem \ref{thm: torsion order and bridge index} in instanton theory:
\begin{enumerate}
	\item The cobordism map as in Theorem \ref{thm: cobordism map, intro}.
	\item The tube attachment lemma as in Proposition \ref{prop: tube attachment lemma}: If $\Sigma'$ is obtained from $\Sigma$ by attaching a tube to $\Sigma$, then (with suitable choices of auxiliary data), we have
	$$F_{\Sigma',A'}=U\circ F_{\Sigma,A}:KHI^-(S^3,K_0)\ra KHI^-(S^3,K_1).$$
	\item A connected sum formula for $KHI^-$.
\end{enumerate}

Item (3) was recently established by the first author and Zemke \cite{GZ2023}, while the first two items will be presented herein.

\subsection{Organization.}The paper is organized as follows. In Section \ref{sec: preliminary} we prove Theorem \ref{thm: main} assuming some technical results that we establish in latter sections. In particular, we construct maps associated to knot cobordisms in Section \ref{sec: cobordism maps} and prove the tube attachment lemma in Section \ref{sec: tube attachment}. 

\subsection{Acknowledgements.} The authors are indebted to Ian Zemke for several helpful discussions about the definition of the knot cobordism map in instanton Floer theory. The idea of the cobordism map developed when the first author was working with Ian Zemke on their paper \cite{GZ2023}. The authors thank John Baldwin, Maggie Miller, Steven Sivek, and Fan Ye for helpful discussions. The first author also thanks the Max Planck Institute for Mathematics for hosting the first author for the bulk of this work.

\section{Preliminaries}\label{sec: preliminary}
\subsection{The minus version of instanton knot homology}
In this section, we review the minus version of instanton knot homology constructed in \cite{li2019direct}. Suppose $Y$ is a connected, closed, oriented $3$-manifold, and let $K\subset Y$ is a null-homologous knot. Let $(\lambda,\mu)$ represent the Seifert framing on $\partial(Y-N(K))$, where $\lambda$ is the longitude and $\mu$ is the meridian. Let $\Gamma_n$ be the suture on $\partial (Y-N(K))$ consisting of two curves of class $\pm(\lambda+n\mu)$. $(Y-N(K),\Gamma_{n})$ forms a balanced sutured manifold as in \cite{juhasz2006holomorphic}, and Kronheimer and Mrowka \cite{kronheimer2010knots} constructed the sutured instanton Floer homology of the pair, denoted by $SHI(Y-N(K),\Gamma_n)$.

A Seifert surface of $K$ induces a grading on $SHI(Y-N(K),\Gamma_n)$. This construction was pioneered by the work of Kronheimer and Mrowka in \cite{kronheimer2010instanton} and Baldwin and Sivek in \cite{baldwin2018khovanov}, and further developed by the second author and his collaborators in \cite{li2019direct,GL2023decompose,LY2022integral1}.

Assuming we fix a Seifert surface $S$ and choose the sutures $\Gamma_n$ to have $2n$ intersections with $S$ for all $n$. In \cite{li2019direct}, if $n$ is odd, then $S$ can be used to construct a $\mathbb{Z}$-grading, while if $n$ is even, one needs to isotope $S$ in one of two ways (called stabilizations in \cite{li2019direct}) to construct a $\mathbb{Z}$-grading. The gradings associated with the two stabilizations always differ by $1$. Therefore, in \cite{LY2022integral1}, the authors avoid choosing one stabilization by defining a $(\mathbb{Z}+\frac{1}{2})$-grading for $S$ when $n$ is even. We will denote the graded part as:
$$SHI(Y-N(K),\Ga_{n},i)~{\rm for~}i\in\intg~{\rm or}~\intg+\frac{1}{2}.$$

As in \cite{baldwin2016contact}, we have two bypass maps
$$\psi^{n}_{\pm,n+1}: SHI(Y-N(K),\Gamma_n)\ra SHI(Y-N(K),\Gamma_{n+1})$$
The original construction involved the orientation reversal of $3$-manifolds and sutures, but we can initiate the construction on $-Y$ to avoid the sign complexity. The two bypass maps are symmetric in some sense, and their grading-shifting behavior was understood in \cite{li2019direct}. In this paper, we adopt the convention that:
$${\rm deg}(\psi^n_{\pm,n+1})=\mp\frac{1}{2}.$$

To simplify the notations, we write
$$\mathbf{\Ga^K_{n}}=SHI(Y-N(K),\Ga_n),~{\rm and~}(\mathbf{\Ga^K_{n}},i)=SHI(Y-N(K),\Ga_n,i).$$
Also we define a grading shift
$$(\mathbf{\Ga^K_{n}}[k],i)=(\mathbf{\Ga^K_{n}},i-k).$$
Under this notation, we know that the map
$$\psi^{n}_{-,n+1}:\mathbf{\Ga^K_{n}}[\frac{1-n}{2}]\to \mathbf{\Ga^K_{n+1}}[\frac{-n}{2}]$$
is grading-preserving and $\psi^{n}_{+,n+1}$ drops the grading by $-1$.
 
Additionally, we can attach a contact $2$-handle (cf. \cite{baldwin2016contact}, and again, we omit the negative orientation of $3$-manifolds) to $(Y-N(K),\Gamma_n)$ along the meridian of the knot, thereby obtaining a map:
$$F_{n}:\mathbf{\Ga^K_{n}}\ra I^{\sharp}(Y).$$

%-------gpt------

\bdefn[\cite{li2019direct,GLW2019}]\label{defn: minus version}
We define $KHI^{-}(Y,K)$ to be the direct limit of the sequence
$$\psi^{n}_{-,n+1}: \mathbf{\Ga^K_{n}}[\frac{1-n}{2}]\to \mathbf{\Ga^K_{n+1}}[\frac{-n}{2}].$$
It inherits a $\intg$-grading. The collection of maps 
$$\psi^{n}_{+,n+1}: \mathbf{\Ga^K_{n}}[\frac{1-n}{2}]\to \mathbf{\Ga^K_{n+1}}[\frac{-n}{2}]$$
induce a map
$$U: KHI^{-}(Y,K,i)\ra KHI^{-}(Y,K,i-1).$$
The maps
$$F_{n}:\mathbf{\Ga^K_{n}}[\frac{1-n}{2}]\ra I^{\sharp}(Y)$$
induce a map
$$F^-_K:KHI^{-}(Y,K)\ra I^{\sharp}(Y).$$
\edefn

\blem[\cite{li2019direct,GLW2019}]\label{lem: minus version for the unknot}
For the unknot $\mathbb{U}_1\subset Y$, we have an isomorphism (of $\mathbb{C}[U]$-modules)
$$KHI^{-}(Y,\mathbb{U}_1)\cong I^{\sharp}(Y)\mathop{\otimes}_{\mathbb{C}}\mathbb{C}[U],$$
where $I^{\#}(Y)$ is supported at grading $0$ and the element $1\in\mathbb{C}[U]$ lies in grading $0$ as well. Furthermore, the map $F^-_{\mathbb{U}_1}$ restricts to an isomorphism on $KHI^{-}(S^3,\mathbb{U}_1,i)$ for all $i\in\intg_{\leq0}$.
\elem

\bdefn
For a null homologous knot $K\subset S^3$, we define
$${\rm ord}_U(K)=\max\{k~|~\exists~x\in KHI^{-}(S^3,K)~s.t.~U^k\cdot x\neq 0,~U^{k+1}\cdot x=0.\}$$
\edefn

\blem\label{lem: odd torsion order}
Suppose $K\subset S^3$ is an alternating knot. Then ${\rm ord}_U(K)$ is odd. 
\elem
\bpf
Suppose $K$ is an alternating knot. In \cite{LY2021large}, the author constructed a set of differentials ${d_k}$ on Kronheimer and Mrowka's version of instanton knot homology $KHI(S^3,K)$. In \cite{LY2022integral2}, it is shown that:
$${\rm ord}_{U}(K)=\max\{k~|~d_{k}\neq0\}.$$
In \cite{LY2021large}, it is demonstrated that $d_k$ has degree $-1$ with respect to the $\mathbb{Z}_2$-grading and degree $-k$ with respect to the Alexander grading on $KHI(S^3,K)$. Furthermore, the $KHI$ of alternating knots is computed through the work of Kronheimer and Mrowka in \cite{kronheimer2010instanton,kronheimer2011khovanov}. Specifically, $KHI(S^3,K,i)$ and $KHI(S^3,K,j)$ are supported in the same $\mathbb{Z}_2$-grading if $i-j$ is even. Hence, $d_k=0$ if $k$ is even, and we are done.
\epf

\subsection{A map for the connected sum of minus version}\label{subsec: connected sum}
\quad
In this section we review a map related to the connected sum of minus version constructed in \cite{GLW2019}. Later it will be used in the construction of the cobordism map. Suppose we have two null-homologous knots $K_1\subset Y_1$ and $K_2\subset Y_2$. In \cite{GLW2019}, the authors constructed a map
$$C_h: \mathbf{\Ga^{K_1}_n}\otimes \mathbf{\Ga^{K_2}_n}\to \mathbf{\Ga^{K_1 \# K_2}_{m+n}}$$

\begin{figure}[ht]
	\begin{overpic}[width=0.8\textwidth]{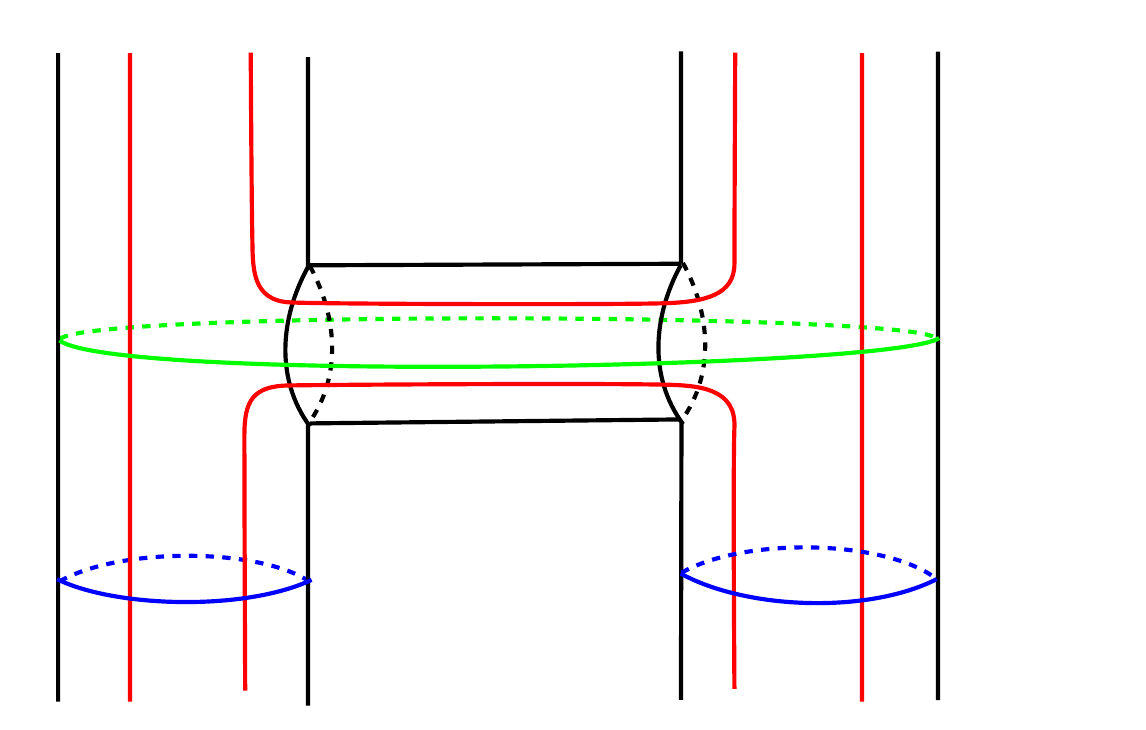}
		\put(42,23){$h^1_1$}
		\put(83,33){$h^2$}
		\put(28,12){$\mu_1$}
		\put(56,12){$\mu_2$}
		\put(8,-2){$Y_1-N(K_1)$}
		\put(64,-2){$Y_2-N(K_2)$}
	\end{overpic}
	\vspace{0.1in}
	\caption{The map $C_h$ is derived through a sequence of operations, initially involving a contact $1$-handle attachment, denoted as $h^1_1$, followed by a contact $2$-handle attachment, represented as $h^2$. The curve guiding the attachment of $h^2$ is depicted as the green curve in the figure provided.}\label{fig: C_h}
\end{figure}
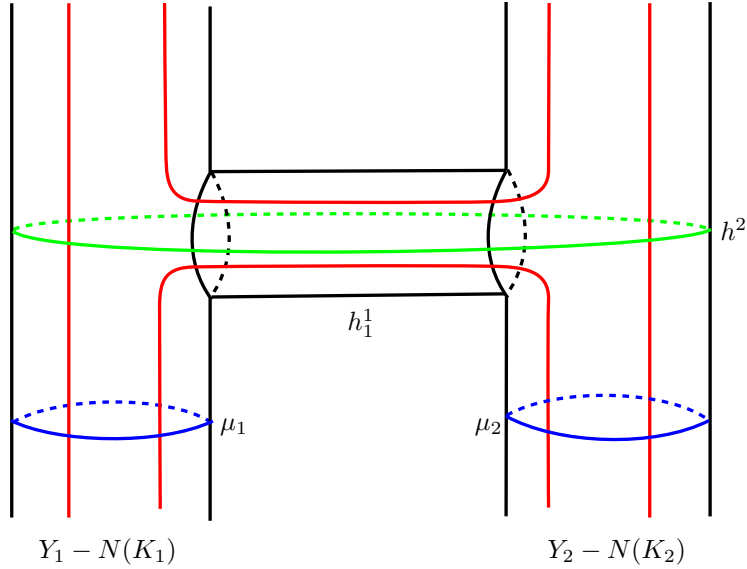

Through a contact $1$-handle attachment followed by a contact $2$-handle attachment, as depicted in Figure \ref{fig: C_h}, the authors also described the grading-shifting behavior of the map $C_h$. Additionally, they prove that the map commutes with the bypass maps on $K_1$ and $K_2$. Essentially, the bypass maps are also obtained from contact handle attachments and are located in regions disjoint from those for $C_h$. Consequently, $C_h$ induces the following map between the minus versions:
\begin{equation}\label{eq: rho map}
	\rho: KHI^-(Y_1,K_1)\otimes KHI^{-}(Y_2,K_2)\to KHI^{-}(Y_1\# Y_2, K_1\# K_2).
\end{equation}

Furthermore, the grading shifting behavior and the commutativity mentioned above also lead to the following two properties of the map $\rho$.

\begin{lem}\label{lem: degree of rho}
	The map $\rho$ is homogeneous with respect to the Alexander grading. In particular,
	$$\rho\bigg(KHI^-(Y_1,K_1,i)\otimes KHI^-(Y_2,K_2,j)\bigg)\subset KHI^-(Y_1\# Y_2, K_1\# K_2,i+j).$$
\end{lem}

\begin{lem}\label{lem: rho intertwines with U}
	The map $\rho$ intertwines with the $U$-action. More precisely, we have
	$$\rho\bigg((U\cdot x)\otimes y\bigg)=\rho\bigg(x\otimes(U\cdot y)\bigg)=U\cdot \rho(x\otimes y).$$
\end{lem}

With the help of the above two lemmas, we can derive the following.
\begin{lem}\label{lem: rho is an iso}
	Suppose $Y_2=S^3$ and $K_2=\mathbb{U}_1$ is the unknot. Let $\nu_0\in KHI^-(Z_0=S^3,B_0=\mathbb{U}_1,0)$ be the generator. Then the map
	$$\iota: KHI^{-}(Y_1,K_1)\to KHI^-(Y_1\cong Y_1\# S^3,K_1\cong K_1\# \mathbb{U}_1)$$
	defined as
	$$\iota(x)=\rho(x\otimes \nu_0)$$
	coincides with the identity map. 
\end{lem}
\bpf
In \cite{li2019direct}, the author demonstrated that the direct system in Definition \ref{defn: minus version} stabilizes when $n$ is large. Transitioning to sutured instanton Floer homology, to validate Lemma \ref{lem: rho is an iso}, it is sufficient to establish the following: let $\nu^m_0 \in (\mathbf{\Gamma}^{\mathbb{U}_1}_m,\frac{m-1}{2})$ be the generator of the highest non-vanishing grading. Then, the map:
$$\iota_{h,n+m}^n: \mathbf{\Ga^{K_1}_n}\to \mathbf{\Ga^{K_1}_{n+m}}\cong \mathbf{\Ga^{K_1\# \mathbb{U}_1}_{m+n}}$$
defined as
$$\iota_{h,n+m}^n(x)=C_h(x\otimes \nu'_0)$$
coincide with the composition of bypass maps:
$$\iota_{h,n+m}^n=\psi^{n+m-1}_{-,m+n}\circ\dots \circ\psi^{n}_{-,n+1}.$$
The sutured instanton Floer homology of the sutured solid tori and the bypass maps on them have been computed in \cite{li2019direct}. In particular, we know that
$$\mathbf{\Ga_{1}^{\mathbb{U}_1}}\cong\mathbb{C},$$
and let $\nu_0^1\in \mathbf{\Ga_{1}^{\mathbb{U}_1}}$ be the generator, then we have
$$\nu_0^m=\psi^{m-1}_{-,m}\circ\dots \circ \psi^{1}_{-,2}(\nu_0^1).$$
Since $C_h$ commutes with the bypass maps, it suffices to prove
$$\iota_{h,n+1}^n=\psi^n_{-,n+1}:\mathbf{\Ga^{K_1}_n}\to \mathbf{\Ga^{K_1}_{n+m}}.$$

%-----GPT-----

To do this, note that $(-(S^3-N(K)),-\Ga_1)$ admits a unique tight contact structure whose contact element is
$$\nu_0^1\in \mathbf{\Ga_{1}^{\mathbb{U}_1}}\cong \mathbb{C}.$$

Note that $(-(S^3-N(K)),-\Gamma_1)$ can be decomposed into a contact $0$-handle and a contact $1$-handle. Additionally, $C_h$ inherently is the map associated with a contact $1$-handle and a contact $2$-handle. In total, we have four contact handles; however, the $0$-handle and one of the $1$-handles form a canceling pair, while the remaining $1$-handle and $2$-handle coincide with the construction of the bypass map as outlined in \cite{baldwin2016contact}. Refer to the accompanying figure for a visual representation. Consequently, by applying the functoriality of the contact gluing map presented in \cite{li2018gluing}, we conclude the following.
$$\iota_{h,n+1}^n=\psi^n_{-,n+1}:\mathbf{\Ga^{K_1}_n}\to \mathbf{\Ga^{K_1}_{n+m}}.$$
\epf

\begin{figure}[ht]
	\begin{overpic}[width=0.8\textwidth]{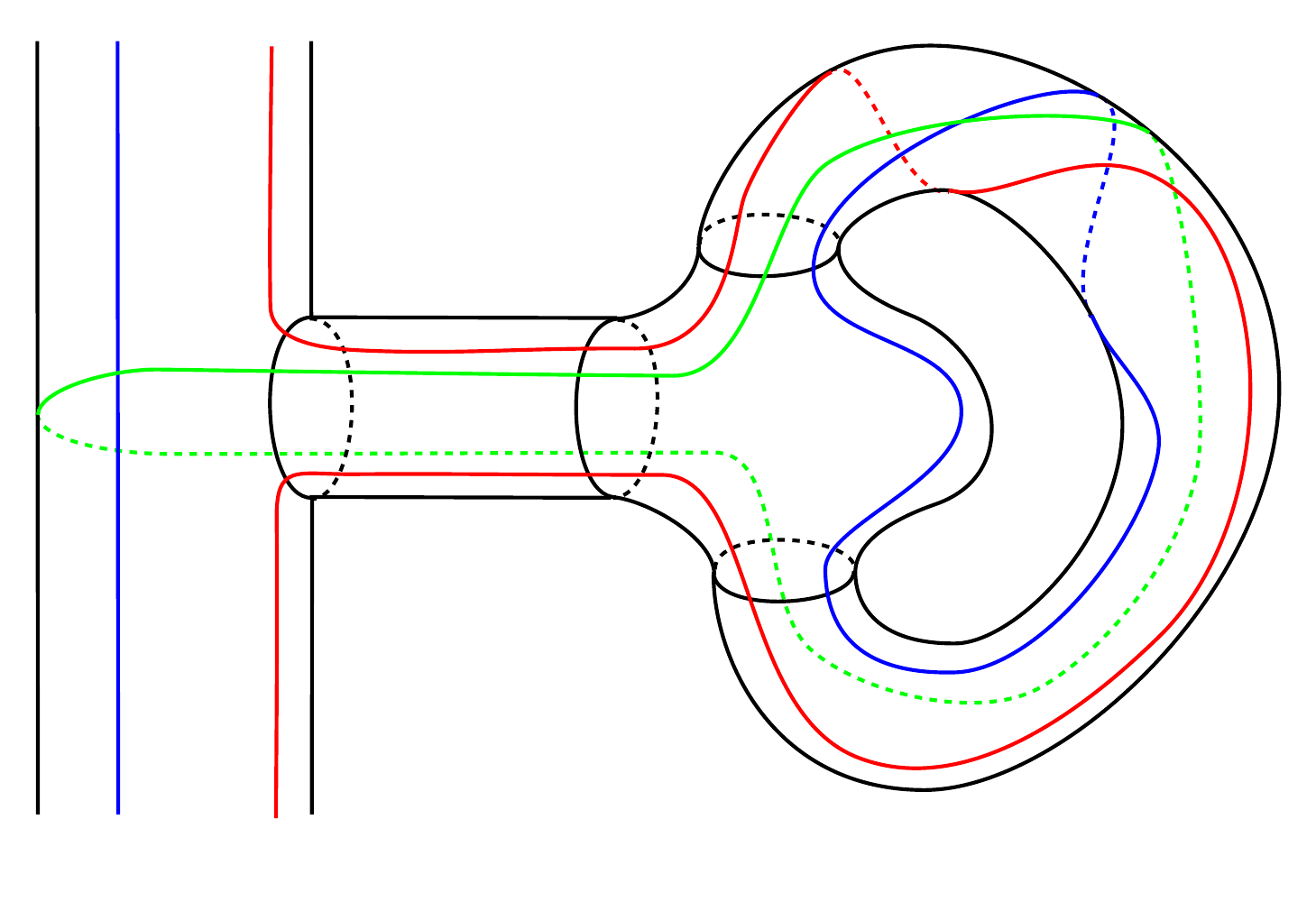}
		\put(35,27.5){$h^1_1$}
		\put(89,59){$h^2$}
		\put(60,38){$h^0$}
		\put(55,11){$h^1_2$}
		\put(8,0){$Y_1-N(K_1)$}
		\put(55,3){$(S^3-N(\mathbb{U}_1))\cong S^1\times D^2$}
	\end{overpic}
	%\vspace{0.1in}
	\caption{The knot complement $S^3-N(K)$ can be identified with a solid torus $S^1\times D^2$, where the meridian of the unknot is identified with the longitude of the solid torus. Then the sutured solid torus can be decomposed into a contact $0$-handle $h^0$ and a contact $1$-handle $h^1_2$. The map $C_h$ is associated with a contact $1$-handle $h^1_1$ and a contact $2$-handle $h^2$, whose attaching curve is the green curve in the picture. The $0$-handle $h^0$ and the $1$-handle $h^1_1$ form a canceling pair, and the $1$-handle $h^1_2$ and the $2$-handle $h^2$ give rise to the bypass map $\psi^{n}_{-,n+1}$.}
\end{figure}
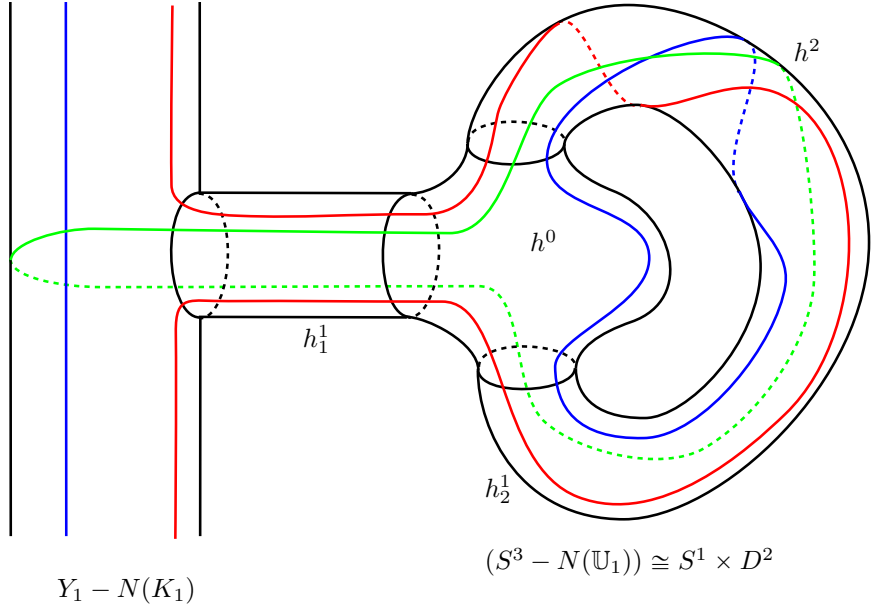

\subsection{Maps associated to special cobordisms}\label{subsec: special cobordism}
In this section we review the definition of special cobordisms the construction of a cobordism map associated to special cobordisms. Suppose we have two null-homologous knots $K_1 \subset Y_1$ and $K_2 \subset Y_2$. Suppose $W$ is a compact, connected, oriented $4$-manifold such that:
$$\partial W=-(Y_1-N(K_1))\mathop{\cup}_{\varphi} (Y_2-N(K_2)),$$
where $\varphi$ identifies the two toroidal boundaries of the knot complements by matching the longitudes and meridians of the two knots. We call such a cobordism a special cobordism, and we can construct a map:
$$F_{W}^-:KHI^{-}(Y_1,K_1)\to KHI^{-}(Y_2,K_2)$$
following \cite{li2018gluing,GLW2019}. We only sketch the ideas here, and details can be found in the papers. We can reinterpret the cobordism $W$ as being obtained from $[0,1] \times (Y_1-N(K_1))$ by attaching a set of 4-dimensional handles along the region ${1}\times{\rm int}(Y_1-N(K_1))$. Hence, this induces a map:
$$F_{W,n}:\mathbf{\Ga^{K_1}_n}\to \mathbf{\Ga^{K_2}_{\varphi(n)}}$$
for all $n$. Note here that $\varphi(n)$ is the slope of the image of $\Gamma_n$ with respect to the framing pushed forward by $\varphi$, and it might be different from the Seifert framing of $K_2$. Note that, as in \cite{baldwin2016contact}, the bypass maps are obtained by gluing contact handles to $\partial (Y_1-N(K_1))$ and are hence disjoint from the handles associated with $W$. As a result, we have the following commutative diagram:
\begin{equation*}
	\xymatrix{
	\mathbf{\Ga^{K_1}_n}\ar[d]^{F_{W,n}}\ar[rr]^{\psi^{n}_{-.n+1}}&&\mathbf{\Ga^{K_1}_{n+1}}\ar[d]^{F_{W,n+1}}\\
	\mathbf{\Ga^{K_2}_n}\ar[rr]^{\psi^{n}_{-.n+1}}&&\mathbf{\Ga^{K_2}_{n+1}}
	}
\end{equation*}
Consequently, by passing to the direct limit, we obtain the map $F_{W}^-$. Essentially, the same argument also substantiates the following:
\blem\label{lem: F^-_W is a module morphism}
The map $F_W^-$ defined as above is a $\mathbb{C}[U]$-module morphism.
\elem
\begin{lem}\label{lem: concordance and 2-handle map}
	Suppose we have two null-homologous knots $K_1\subset Y_1$ and $K_2\subset Y_2$. Suppose $W'$ is a cobordism from $Y_1$ to $Y_2$ obtained from $[0,1]\times Y_1$ by attaching a set of 4-dimensional handles along the region ${1}\times \text{int}(Y_1-N(K_1))$. The same set of handles also defines a special cobordism, $W$, from $Y_1-N(K_1)$ to $Y_2-N(K_2)$. Then we have a commutative diagram:
\begin{equation*}
	\xymatrix{
	KHI^-(Y_1,K_1)\ar[rr]^{F^-_{W}}\ar[d]^{F_{K_1}^-}&&KHI^-(Y_2,K_2)\ar[d]^{F_{K_2}^-}\\
	I^{\#}(Y_1)\ar[rr]^{F_{W'}}&&I^{\#}(Y_2)
	}
\end{equation*}
\end{lem}

\subsection{Borromean knots and the $H_1$-action on instanton Floer homology}\label{subsec: borromean knots}
\begin{figure}[ht]
	\begin{overpic}[width=0.4\textwidth]{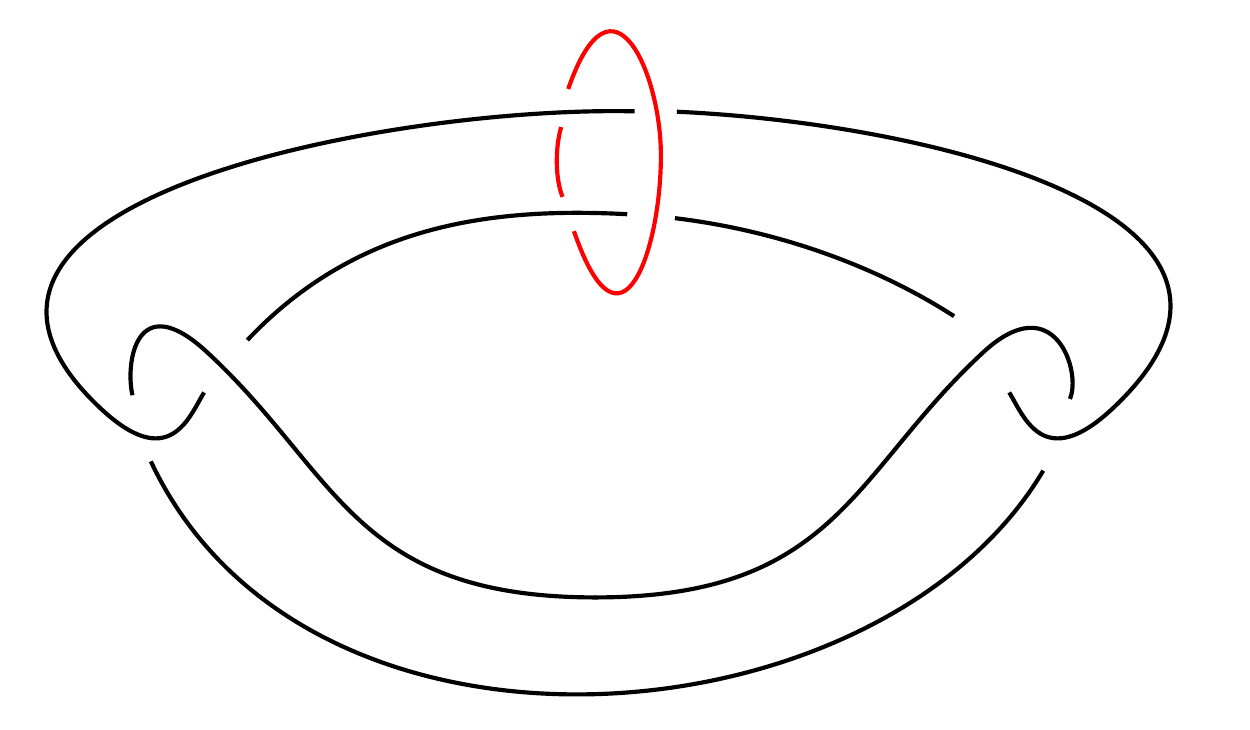}
		%\put(10,19){$\delta$}
		%\put(16,2){$K_+$}
		\put(47,30){\color{red} $B_1$}
		\put(80,46){$0$}
		\put(80,10){$0$}
	\end{overpic}
	\caption{The Borromean knot $B_1$ is inside $Z_2 = S^1 \times S^2 \# S^1 \times S^2$. The two copies of $S^1 \times S^2$ arise from the zero surgeries on the two (black) components of the Borromean link.\label{fig: borromean}}
\end{figure}

In this section we review some basic facts related to the Borromean knots, which will later be used in the construction of cobordism maps for $KHI^-$. For any $n \in \intg_{\geq 0}$, we denote $Z_n = \#^n(S^1 \times S^2)$. Let $B_1 \subset Z_2$ be the Borromean knot, as depicted in \ref{fig: borromean}. Let $B_n = \#^n B_1 \subset Z_{2n}$, representing the connected sum of $2n$ copies of $B_1$. As indicated in \cite[Section 9]{ozsvath2004holomorphicknot}, the knot complement $Z_2 - N(B_1)$ is diffeomorphic to $S^1 \times (T^2 - D^2)$. On the boundary, the meridian of $B_1$ aligns with $S^1$, and the Seifert longitude of $B_1$ corresponds to $\partial D^2$. We make the following observation.

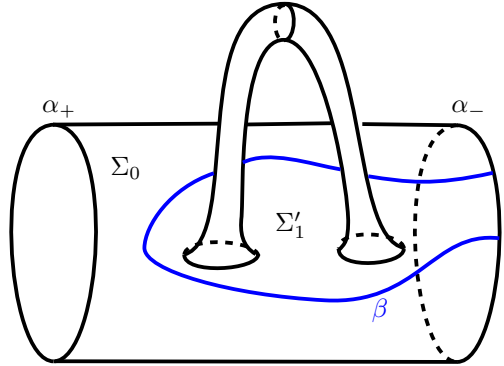
\begin{figure}[ht]
	\begin{overpic}[width=0.5\textwidth]{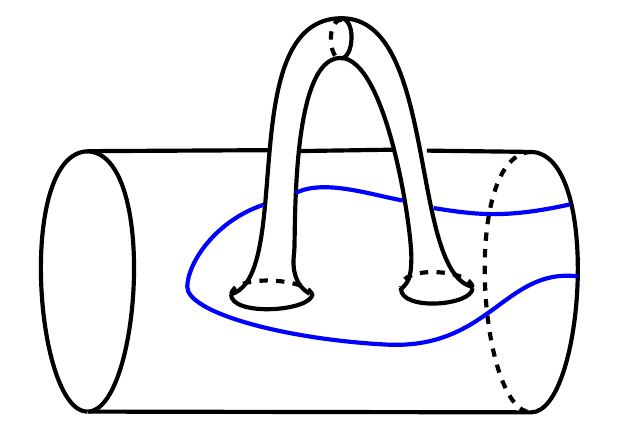}
		%\put(10,19){$\delta$}
		%\put(16,2){$K_+$}
		\put(53,27){$\Sigma'_1$}
		\put(24,37){$\Sigma_0$}
		\put(84,48){$\al_-$}
		\put(12,48){$\al_+$}
		\put(70,12){\color{blue}$\beta$}
	\end{overpic}
	\caption{The surface $\Sigma_1$ and the curve $\beta$.\label{fig: Sigma_1}}
\end{figure}

\blem\label{lem: connected sum with B_g}
Suppose $K \subset Y$ is a null-homologous knot and $(\lambda, \mu)$ constitutes the Seifert framing on the boundary of the knot complement. Let $\Sigma_g$ be a compact, connected, oriented surface of genus $g$ with two boundary components $\alpha_+$ and $\alpha_-$. Additionally, assume we have a diffeomorphism
$$f: \partial (Y-N(K))\ra S^1\times\al_+{~s.t.}~f(\lambda)=\al_+,~{\rm and~}f(\mu)=S^1.$$
Then, we have a diffeomorphism
$$\bigg((Y-N(K))\mathop{\cup}_f(S^1\times\Sigma_g)\bigg)\cong\bigg(Y\#Z_{2g}-N(K\#B_g)\bigg).$$
Furthermore, if $(\lambda_{\sharp},\mu_{\sharp})$ is the Seifert framing of the connected sum knot $K\# B_g$, then on the boundary the above diffeomorphism identifies $S^1$ with $\mu_{\sharp}$ and $\al_-$ with $\lambda_{\sharp}$.
\elem
\bpf
We first prove the case $g=1$. As in Figure \ref{fig: Sigma_1}, we pick a properly embedded arc $\beta\subset \Sigma_1$ dividing it into two parts
$$\Sigma_1=\Sigma_0\mathop{\cup}_{\beta}\Sigma_1',$$
where $\Sigma_0$ is an annulus and $\al_+\subset \partial \Sigma_0$, and $\Sigma_1'$ is a once-punctured torus. As a result, we know 
$$\bigg((Y-N(K))\mathop{\cup}_f(S^1\times\Sigma_1)\bigg)=\bigg((Y-N(K))\mathop{\cup}_f(S^1\times\Sigma_0)\bigg)\mathop{\cup}_{S^1\times \beta}\bigg(S^1\times\Sigma_1'\bigg).$$
Clearly, $(Y-N(K))\mathop{\cup}_f(S^1\times\Sigma_0)$ is diffeomorphic to the knot complement $Y-N(K)$; $S^1\times\Sigma_1'$ is diffeomorphic to the knot complement $Z_2-N(B_1)$; and $S^1\times\beta$ is an annular neighborhood of the meridian on the two boundaries. Hence we are done for $g=1$.

For a general $g$, we can decompose $\Sigma_g$ into $g$ copies of $\Sigma_1$ and use induction to complete the proof.
\epf

%For a closed $3$-manifold $Y$ and a homology class $\al\in H_*(Y)$, there is a standard construction of $\mu$-action
%$$\mu(\al):I^{\sharp}(Y)\ra I^{\sharp}(Y).$$

The instanton Floer homology of connected sums of $S^1\times S^2$ have been computed in {\cite[Section 7.8]{scaduto2015instanton}}. In particular, we have the following.

\begin{prop}\label{prop: H_1-action on S1 times S2}
	Suppose $Z$ is a closed manifold that is diffeomorphic to connected sums of $n$ copies of $S^1\times S^2$. Then there is a well-defined isomorphism
	$$\phi_n: \Lambda^* H_1(Z;\mathbb{C})\xra{\cong} I^{\sharp}(Z).$$
\end{prop}

The minus version $KHI^-$ of the Borromean knot $B_g\subset Z_{2g}$ has been computed by the second author and Ye in \cite{LY2022integral2}. In particular, we have the following.

\bprop\label{prop: KHI- for Borromean knots}
Suppose $B_g\subset Z_{2g}$ is defined as above. Then for any grading $j$, the map
$$F_{B_g,j}^-:KHI^{-}(Z_{2g},B_{g},j)\to I^{\sharp}(Z_{2g})$$
is injective with image
$$F_{B_g,j}^-\bigg(KHI^{-}(Z_{2g},B_{g},j)\bigg)=\phi_{2g}\bigg(\bigoplus_{0\leq k\leq g-j}\Lambda^k H_1(Z_{2g};\mathbb{C})\bigg)$$
Furthermore, under these identifications, the $U$-action 
$$U: KHI^{-}(Z_{2g},B_{g},j)\ra KHI^{-}(Z_{2g},B_{g},j-1)$$
coincide with the inclusion
$$\phi_{2g}\bigg(\bigoplus_{0\leq k\leq g-j}\Lambda^k H_1(Z_{2g};\mathbb{C})\bigg)\hookrightarrow \phi_{2g}\bigg(\bigoplus_{0\leq k\leq g-j+1}\Lambda^k H_1(Z_{2g};\mathbb{C})\bigg)$$
\eprop 

Now with the identification from above Propositions, we identify $I^{\sharp}(Z_{2g})$ with $\Lambda^*H_1(Z_{2g};\mathbb{C})$ and omit the isomorphism $\phi_{2g}$. We can take a set of generators $\theta_1,...,\theta_{2g}$ of $H_1(Z_{2g};\mathbb{C})$ and view them as elements of $I^{\sharp}(Z_{2g})$. Note that for any element $\al\in H_1(Z_{2g})$, there is a $\mu(\al)$-action
$$\mu(\al):I^{\#}(Z_{2g})\to I^{\#}(Z_{2g}).$$
As in \cite{LY2022integral2}, we can describe the $H_1(Z_{2g})$-action on $I^{\#}(Z_{2g})$ as follows: we can choose a second set of basis $\iota_1,...,\iota_{g}$ of $H_1(Z_{2g};\mathbb{C})$, and treat the $H_1(Z_{2g})$-action on $I^{\sharp}(Z_{2g})$ as the contraction:
$$\iota_{i}(\theta_j)=\begin{cases}
	1&i=j\\
	0&i\neq j
\end{cases},~{\rm and}~\iota_i(\al\wedge\beta)=\iota_i(\al)\wedge\be+(-1)^{{\rm deg}(\al)}\al\wedge\iota_i(\beta).$$

%-----GPT-----

\bdefn\label{defn: the element nu_g}
For the knot $B_g\subset Z_{2g}$, define
$$\nu_g=F_{B_{g},-g}^{-1}(\theta_1\wedge...\wedge\theta_{2g})\in KHI^{-}(Z_{2g},B_g,-g).$$
\edefn

\brem
Note, as mentioned in Remark \ref{rem: scalar ambiguity}, that sutured instanton Floer homology is well-defined only up to a scalar. For different choices of basis elements $\theta_1$,..., $\theta_{g}$, their total wedge product differs by a scalar ($\Lambda^gH_1(Z_{2g};\mathbb{C})$ is $1$-dimensional). If we want to make an explicit choice of the basis elements $\theta_i$ as well as $\nu_g$, we can refer to {\cite[Section 7.8]{scaduto2015instanton}}, where all such elements in $I^{\#}(Z_{2n})$ are chosen as the relative invariants of suitably chosen four manifolds bounded by $Z_{2n}$.
\erem

Note that for any $g_1,g_2\in\intg_{\geq0}$, we have
$$(Z_{2(g_1+g_2)},B_{g_1+g_2})=(Z_{2g_1},B_{g_1})\#(Z_{2g_2},B_{g_2}),$$
and hence we have a map as in Section \ref{subsec: connected sum}
$$\rho:KHI^{-}(Z_{2g_1},B_{g_1})\otimes KHI^{-}(Z_{2g_2},B_{g_2})\to KHI^{-}(Z_{2(g_1+g_2)},B_{g_1+g_2}).$$
\begin{lem}\label{lem: tensor product of nu}
	We have the following identity
	$$\rho(\nu_{g_1}\otimes \nu_{g_2})=\nu_{g_1+g_2}.$$
\end{lem}
\bpf
In \cite{GLW2019}, the authors proved the following commutative diagram
\begin{equation*}
	\xymatrix{
	KHI^{-}(Z_{2g_1},B_{g_1})\otimes KHI^{-}(Z_{2g_2},B_{g_2})\ar[rr]^{\rho}\ar[d]^{F_{B_{g_1}}^-\otimes F^-_{B_{g_2}}}&&KHI^{-}(Z_{2(g_1+g_2)},B_{g_1+g_2})\ar[d]^{F^-_{B_{g_1+g_2}}}\\
	I^{\#}(Z_{2g_1})\otimes I^{\#}(Z_{2g_2})\ar[rr]^{F_{\#}}&&I^{\#}(Z_{2g_1+2g_2})
	}
\end{equation*}
Here, the map
$$F_{\#}:I^{\#}(Z_{2g_1})\otimes I^{\#}(Z_{2g_2})\to I^{\#}(Z_{2g_1+2g_2})
	$$
comes from taking the connected sum of two manifolds. The fact
$$F_{\#}\bigg((\theta_1\wedge...\wedge\theta_{2g_1})\otimes (\theta'_1\wedge...\wedge\theta'_{2g_2})\bigg)=\theta_1\wedge...\wedge\theta_{2g_1}\wedge\theta'_1\wedge...\wedge\theta'_{2g_2}$$
follows from \cite[Section 7.8]{Scaduto2018}. (Essentially, the elements in $I^{\#}(Z_{n})$ for ($n=g_1,g_2,g_1+g_2$) are the relative invariants from boundary connected sums of $n$ copies of $S^1\times D^3$, and $F_{\#}$ maps the tensor product of such elements to such elements.)
\epf

Now recall $Z_{n}$ is the connected sum of $n$ copies of $S^1\times S^2$ (and $Z_0=S^3$). $Z_n$ can be obtained from $Z_{n-1}$ by a $0$-surgery along an unknot. As a result, we can form a cobordism $W^n_{n-1}$ from $Z_n$ to $Z_{n-1}$.  For $n>m\geq 0$ write 
$$W^{n}_{m}=W^{m+1}_{m}\cup...\cup W^{n}_{n-1}.$$
Let
$$G^n_m:I^{\sharp}(Z_n)\ra I^{\sharp}(Z_m)$$
be the associated cobordism map.
We have the following.
\blem\label{lem: G is surjective}
For any $n>m\geq 0$, the map $G^n_m$ is surjective.
\elem
\bpf
For any $k$, the map $G^{k}_{k-1}$ fits into an exact sequence
\begin{equation*}
	\xymatrix{
	I^{\sharp}(Z_k)\ar[rr]&&I^{\sharp}(Z_k)\ar[dl]\\
	&I^{\sharp}(Z_{k+1})\ar[ul]^{G^{k+1}_{k}}&
	}
\end{equation*}
Hence the surjectivity follows from Proposition \ref{prop: H_1-action on S1 times S2}.
\epf

Note that when we obtain $Y_{2g}$ from $S^3$ by performing 0-surgeries along the unlink $\mathbb{U}_{2g}$, there exists an unknot $\mathbb{U}_1 \subset S^3$, disjoint from $\mathbb{U}_{2g}$, that becomes the knot $B_g \subset Z_{2g}$. Consequently, the $0$-surgeries give rise to a special cobordism $X^n_m$ from $Z_{2n}-N(B_n)$ to $Z_{2m}-N(B_m)$. Here, $B_0 = \mathbb{U}_1 \subset Z_0 = S^3$. As in Section \ref{subsec: special cobordism}, we then obtain a map.
$$F^n_m=F^-_{X^n_m}:KHI^{-}(Z_{2n},B_n)\ra KHI^{-}(Z_{2m},B_m).$$
\blem\label{lem: the map F^2n_0}
For the map $F^{2g}_0$ defined as above, we have 
$$F^{2g}_0(\nu_g)=U^{g}\in KHI^{-}(S^3,\mathbb{U}_1,-g).$$
Furthermore, for $k<2g$, we have
$$F^{2g}_0(\theta_{i_1}\wedge...\wedge\theta_{i_k})=0.$$
\elem
\bpf
From Lemma \ref{lem: concordance and 2-handle map}, there is a commutative diagram
\begin{equation*}
	\xymatrix{
	KHI^{-}(Z_{2g},B_g)\ar[rr]^{F^{2g}_0}\ar[d]^{F^-_{B_g}}&&KHI^{-}(S^3,\mathbb{U}_1)\ar[d]^{F^-_{\mathbb{U}_1}}\\
	I^{\sharp}(Z_{2g})\ar[rr]^{G^{2g}_0}&& I^{\sharp}(S^3)
	}
\end{equation*}
Note the map $G^{2g}_0$ is also equivariant with the $H_1$-action. Since $H_1(S^3)=0$, we know for $i=1,...,2g$, $G^{2g}_0\circ\iota_i=0$. As a result, from the above commutative diagram and the fact that $F^-_{\mathbb{U}_1}$ is an isomorphism on $KHI^{-}(S^3,\mathbb{U}_1,j)$ for any $j$ ({\it c.f.} Lemma \ref{lem: minus version for the unknot}), we know that.
$$F^{2g}_0(\theta_{i_1}\wedge...\wedge\theta_{i_k})=0$$
For any $k<2g$. Furthermore, Lemma \ref{lem: G is surjective} implies that $F^{2g}_0(\nu_g)\neq 0$. Note the unknot $\mathbb{U}_{1}$ bounds a genus-$g$ surface in the link complement $S^3-N(\mathbb{U}_{2g})$, which becomes a minimal genus Seifert surface of $B_g\subset Z_{2g}$. We conclude that the map $F^{2g}_0$ preserves grading. Consequently, we conclude that
$$F^{2g}_0(\nu_g)=U^{g}\in KHI^{-}(S^3,\mathbb{U}_1,-g).$$
\epf

\section{Maps associated to knot cobordisms}\label{sec: cobordism maps}
In this section, we describe cobordism maps associated with various types of cobordisms between knots.

\subsection{A relative invariant}\label{subsec: one end} 
Assume that $W$ is a connected, compact, oriented $4$-manifold such that $\partial W = Y$. Let $\Sigma_g$ be a connected, compact, oriented surface properly embedded, so that
$$K=\partial \Sigma_g\cap Y$$
is a null-homologous knot. Then we can associate an element
$$\theta(W,\Sigma_g)\in KHI^-(Y,K)$$
For the pair $(W,\Sigma_g)$, the procedure is as follows. Take $M=S^1\times \Sigma_g$. There is a canonical framing on $\partial M\cong T^2$, where the meridian is the $S^1$-direction and the longitude is the $\partial \Sigma_g$-direction. As in Section \ref{subsec: borromean knots}, if $Z_{2g}$ is the Dehn filling of $M$ along meridian, then $Z_{2g}$ is diffeomorphic to the connected sum of $2g$ copies of $S^1\times S^2$ and the core of the solid torus becomes the Borromean knot $B_{2g}\subset Z_{2g}$. Now as in Proposition \ref{prop: H_1-action on S1 times S2}, there is a well-defined isomorphism
$$\phi_{2g}:\Lambda^*H_1(Z_{2g};\mathbb{C})\to I^{\#}(Z_{2g}).$$
Hence we can pick the element $\nu_{g}\in KHI^-(Z_{2g},B_{g})$ as in Definition \ref{defn: the element nu_g}. As in Section \ref{subsec: special cobordism}, $W-N(\Sigma_g)$ can be viewed as a special cobordism from the Borromean knot complement $Z_{2g}-N(B_g)$ to the knot complement $Y-N(K)$. This leads to a map:
$$F^-_{W-N(\Sigma_g)}:KHI^-(Z_{2g},B_g)\to KHI^-(Y,K).$$
We can then define
\begin{equation}\label{eq: relative invariant}
	\theta(W,\Sigma_g)=F^-_{W-N(\Sigma_g)}(\nu_g),
\end{equation}

\subsection{Cobordism between two knots}\label{subsec: two ends} Assume that for $i=0,1$, we have a null-homologous knot $K_i\subset Y_i$, and $W$ is a connected oriented $4$-manifold so that
$$\partial W=-Y_0\cup Y_1.$$

\begin{figure}[ht]
	\begin{overpic}[width=0.8\textwidth]{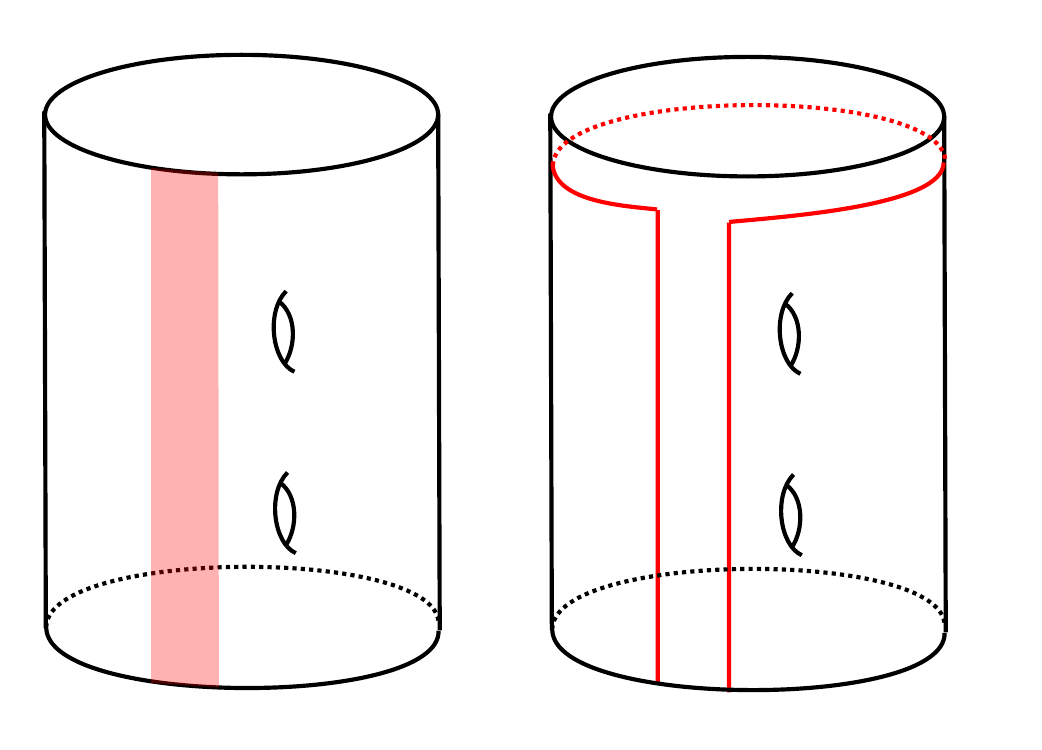}
		\put(28,57){$K_0$}
		\put(28,7.5){$K_1$}
		\put(16.5,30){$\mathcal{D}$}
		\put(80,49){$\beta$}
		\put(80,30){$\Sigma_g'$}
		\put(65,30){$\Sigma_0$}
	\end{overpic}
	\caption{The shaded region is the decoration $\mathcal{D}$, and the red curve is $\beta$.}\label{fig: decoration_1}
\end{figure}
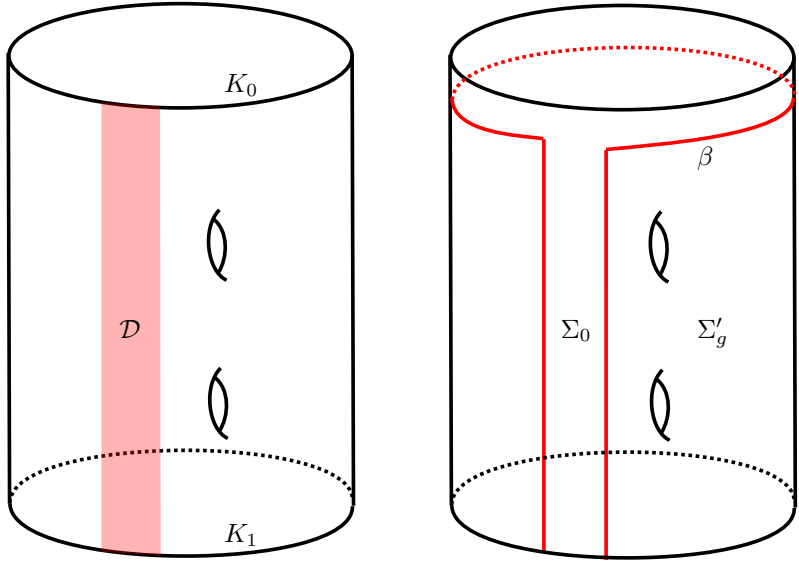

\bdefn\label{defn: Specially decorated}
A \textit{specially decorated} knot cobordism $(W, \Sigma_g,\mathcal{D}): (Y_0, K_0) \to (Y_1, K_1)$ between two knots consists of a pair $(W, \Sigma_g)$ where $\Sigma \subset W$ is a properly embedded surface with $\Sigma\cap Y_i=K_i$, and a codimension-$0$ sub-manifold $\mathcal{D}\subset \Sigma_g$ equipped with an embedding
$$\eta:[0,1]\times[0,1]\to \mathcal{D}$$
with
$$\mathcal{D}\cap K_0=\eta(\{0\}\times[0,1])~{\rm and~}\mathcal{D}\cap K_1=\eta(\{1\}\times[0,1]).$$
\edefn
See Figure \ref{fig: decoration_1}.

We define a map
$$F_{W,\Sigma_g,\mathcal{D}}: KHI^-(Y_0,K_0)\to KHI^-(Y_1,K_1)$$
as follows. Let 
$$\beta'=\eta([0,1]\times\{0,1\})\cup (K_0-\eta(\{0\}\times [0,1]))$$
and $\beta$ be obtained from $\beta'$ by smoothing the corner and being pushed into the interior of $\Sigma_{g}-\mathcal{D}$. See Figure \ref{fig: decoration_1} The curve $\beta$ cuts the surface $\Sigma_{g}$ into two parts:
$$\Sigma_g=\Sigma_0\mathop{\cup}_{\beta}\Sigma'_g$$
where $\Sigma_0$ is an annulus with a canonical identification $\Sigma_0=[0,1]\times K_0$ and $\Sigma'_g$ is its complement. Note that, as in Lemma \ref{lem: connected sum with B_g},
$$S^1\times \Sigma'_g\cong Z_{2g}-N(B_g)$$
and $\bigg((Y_0-N(K_0))\cup (S^1\times \Sigma_g)\bigg)$ can be identified as the complement of the connected sum knot:
$$\bigg((Y_0-N(K_0))\cup (S^1\times \Sigma_g)\bigg)\cong \bigg(Y_0\# Z_{2g}-N(K_0\# B_g)\bigg),$$
where the connected sum is performed along the annular region $S^1\times\beta$. As a result, we have the map
$$\rho:KHI^-(Y_0,K_0)\otimes KHI^-(Z_{2g},B_g)\to KHI^-(Y_0\# Z_{2g},K_0\# B_{2g}).$$
Take $\nu_g\in KHI^{-}(Z_{2g},B_g)$ as in Definition \ref{defn: the element nu_g}. The $4$-manifold $W-N(\Sigma_g)$ can be viewed as a special cobordism from $Y_0\# Z_{2g}-N(K_0\# B_g)$ to $Y_1-N(K_1)$, so it induces a special cobordism and hence a map between the minus versions as in Section \ref{subsec: special cobordism}. We can finally define $F_{W,\Sigma_g,\mathcal{D}}$ as:
\begin{equation}\label{eq: cobordism map from K_0 to K_1}
	F_{W,\Sigma_g,\mathcal{D}}(x)=F^-_{W-N(\Sigma_g)}\circ\rho(x\otimes \nu_g).
\end{equation}

\blem\label{lem: F_W,Sigma,D is a module morphism}
The map $F_{W,\Sigma_g,\mathcal{D}}$ is a $\mathbb{C}[U]$-module morphism.
\elem
\bpf
It follows from Lemma \ref{lem: rho intertwines with U} and Lemma \ref{lem: F^-_W is a module morphism}.
\epf

\blem\label{lem: functoriality of F_W,Sigma,D}
Suppose we have a cobordism $(W,\Sigma_{g},\mathcal{D})$ from $(Y_0,K_0)$ to $(Y_1,K_1)$ and a cobordism $(W',\Sigma'_{g'},\mathcal{D}')$ from $(Y_1,K_1)$ to $(Y_2,K_2)$. Suppose further that $\mathcal{D}\cup \mathcal{D}'$ is also a decoration on $\Sigma_g\cup\Sigma'{g'}$. Then we have the following identity:
$$F_{W\cup W',\Sigma_g\cup\Sigma'_{g'},\mathcal{D}\cup\mathcal{D}'}=F_{W',\Sigma'_{g'},\mathcal{D}'}\circ F_{W,\Sigma,\mathcal{D}}.$$
\elem

\bpf
The proof of this lemma is essentially the same as the argument in \cite{li2018gluing}. From the definition of the cobordism map, we know
$$F_{W',\Sigma'_{g'},\mathcal{D}'}\circ F_{W,\Sigma,\mathcal{D}}(x)=F^-_{W'-N(\Sigma'_{g'})}\circ\rho'\bigg((F_{W-N(\Sigma_g)}^-\circ\rho(x\otimes \nu_g))\otimes \nu_{g'}\bigg),$$
where
$$\rho:KHI^-(Y_0,K_0)\otimes KHI^-(Z_{2g},B_{2g})\to KHI^{-}(Y_0\# Z_{2g}, K_0\# B_{g}),$$
and
$$\rho': KHI^-(Y_1,K_1)\otimes KHI^-(Z_{2g'},B_{2g'})\to KHI^{-}(Y_1\# Z_{2g'}, K_1\# B_{g'}).$$
Note that the handles involved in constructing the map $F_{W-N(\Sigma_g)}^-$ are glued to the region 
$${\rm int}\bigg((Y_0-N(K_0)\cup(S^1\times \Sigma^g)\bigg),$$ while the map $\rho'$ is constructed via gluing contact handles to
$$\partial \bigg(Y_1-N(K_1)\bigg)=\partial \bigg((Y_0-N(K_0)\cup(S^1\times \Sigma^g)\bigg).$$
Hence, essentially the two maps $\rho'$ and $F_{W-N(\Sigma_g)}^-$ commute, and we have
\beq
F_{W',\Sigma'_{g'},\mathcal{D}'}\circ F_{W,\Sigma,\mathcal{D}}(x)&=F^-_{W'-N(\Sigma'_{g'})}\circ\rho'\bigg((F_{W-N(\Sigma_g)}^-\circ\rho(x\otimes \nu_g))\otimes \nu_{g'}\bigg)\\
&=F^-_{W'-N(\Sigma'_{g'})}\circ F_{W-N(\Sigma_g)}^-\circ \rho''\bigg(\rho(x\otimes\nu_{g})\otimes\nu_{g'}\bigg),
\eeq
where we have
$$\rho'':KHI^-(Y_0\# Z_{2g},K_{0}\#B_g)\otimes KHI^-(Z_{2g'},B_{2g'})\to KHI^-(Y_0\# Z_{2g}\# Z_{2g'},K_{0}\#B_g\# B_{g'}).$$
Note that the contact handles associated to $\rho''$ and $\rho$ are attached in disjoint regions so they commute, and we obtain the following equality
\beq
\rho''\bigg(\rho(x\otimes\nu_{g})\otimes\nu_{g'}\bigg)&=\rho^{(3)}\bigg(x\otimes\rho^{(4)}(\nu_{g}\otimes\nu_{g'})\bigg)\\
[{\rm by~Lemma~\ref{lem: tensor product of nu}}]~&=\rho^{(3)}(x\otimes \nu_{g+g'}),
\eeq
where we have
$$\rho^{(3)}:KHI^-(Y_0,K_0)\otimes KHI^{-}(Z_{2(g+g')},B_{g+g'})\to KHI^-(Y_0\# Z_{2(g+g')},K_{0}\#B_{g+g'})$$
and
$$\rho^{(4)}:KHI^-(Z_{2g},B_{g})\otimes KHI^-(Z_{2g'},B_{g'})\to KHI^-(Z_{2(g+g')},B_{g+g'}).$$
Finally, the handles involved in $F_{W'-N(\Sigma'_{g'})}^-$ and $F^-_{W-N(\Sigma_{g})}$ together give rise to handles involved in $F^-_{(W\cup W')-N(\Sigma_g\cup\Sigma'_{g'})}$. As a result, we have
\beq
F_{W',\Sigma'_{g'},\mathcal{D}'}\circ F_{W,\Sigma,\mathcal{D}}(x)&=F^-_{W'-N(\Sigma'_{g'})}\circ F_{W-N(\Sigma_g)}^-\circ \rho''\bigg(\rho(x\otimes\nu_{g})\otimes\nu_{g'}\bigg)\\
&=F^-_{(W\cup W')-N(\Sigma_g\cup\Sigma'_{g'})}\circ \rho^{(3)}(x\otimes \nu_{g+g'})\\
&=F_{W\cup W',\Sigma_g\cup\Sigma'_{g'},\mathcal{D}\cup\mathcal{D}'}(x).
\eeq
\epf

\subsection{Three ends}\label{subsec: three ends} Suppose that for $i=0,1,2$, we have null-homologous knots $K_i\subset Y_i$ such that $K_1\subset Y_1$ is the unknot. Suppose $W$ is a compact, connected, oriented $4$-manifold such that
$$\partial W = -Y_0\cup -Y_1\cup Y_2$$
and $\Sigma_g\subset W$ is a properly embedded surface such that $\Sigma_g\cap Y_i = K_i$ for $i=0,1,2$. We also want a decoration on $\Sigma_g$. Let $\mathcal{D}\subset \Sigma_g$ be a sub-manifold such that there is a diffeomorphism
$$\eta: H\to \mathcal{D},$$
where $H$ is a regular hexagon whose six edges are labeled as $e_1$,..., $e_6$ in order, with the following property: For $i=0,1,2$, we have
$$\mathcal{D}\cap K_i = \eta(e_{2i+1}).$$

Take $\beta' = \eta(e_6)\cup (K_0-\eta(e_1))\cup \eta(e_2\cup e_3\cup e_4)$, and $\delta' = \beta'-\eta(e_3)\cup (K_1-\eta(e_3))$. We smooth and push $\beta'$ and $\delta'$ to obtain properly embedded curves $\beta$ and $\delta$, respectively, so that
$$\Sigma_g=\Sigma_{0}\cup \Sigma'_{0}\cup\Sigma'_{g},$$
where $\Sigma_{0}$ and $\Sigma'_{0}$ are both annuli, $K_0\subset \partial \Sigma_0$, and $K_1\subset \partial \Sigma'_{0}$. This data specifies an identification
$$\bigg((Y_0-N(K_0))\sqcup(Y_1-N(K_1))\cup S^1\times \Sigma_g\bigg)\cong\bigg(Y_0\# Y_1\# Z_{2g}-N(K_0\# K_1\# B_g)\bigg).$$

\begin{figure}[ht]
	\begin{overpic}[width=1\textwidth]{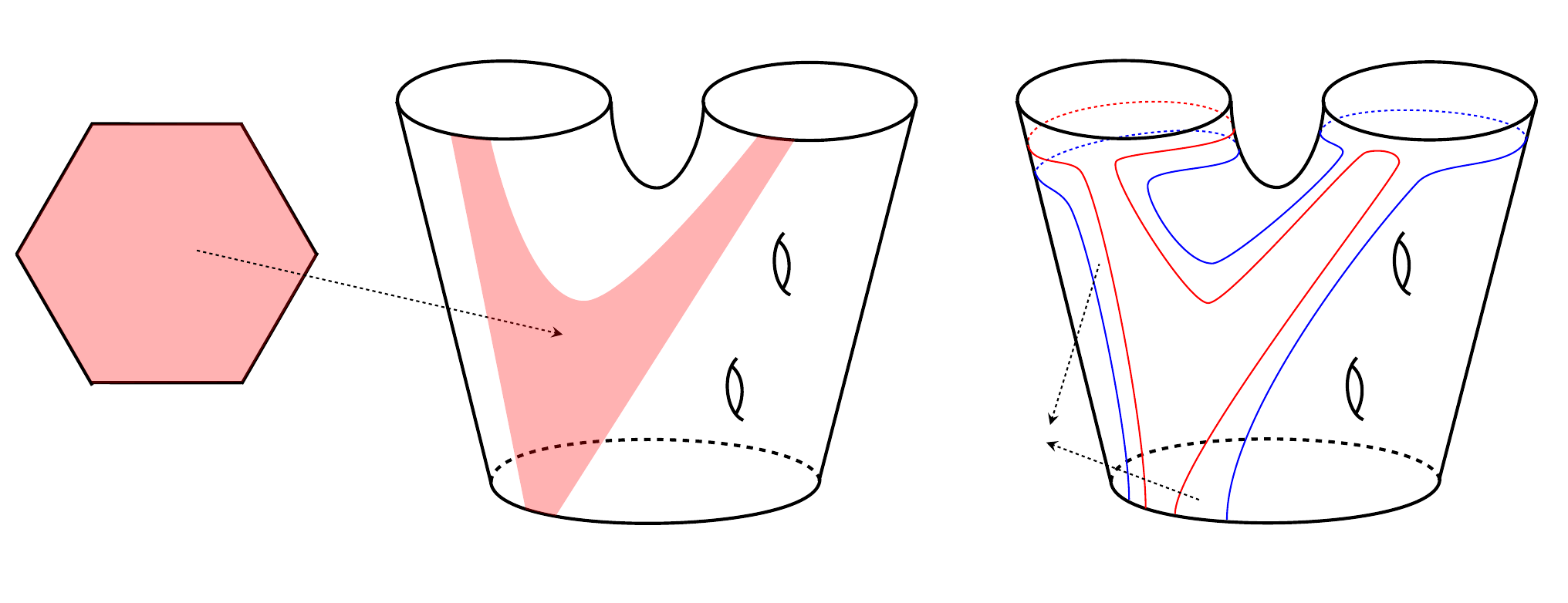}
		\put(9.5,32){$e_1$}
		\put(18,16){$e_3$}
		\put(1,16){$e_5$}
		\put(31,36){$K_0$}
		\put(51,36){$K_1$}
		\put(41,3){$K_2$}
		\put(9.5,22){$H$}
		\put(36.5,16){$\mathcal{D}$}
		\put(80,49){$\beta$}
		\put(89,15){$\Sigma_g'$}
		\put(63.5,10.5){$\Sigma_0'$}
		\put(75,15){$\Sigma_0$}
		\put(80,20){\color{red}$\beta$}
		\put(95,26){\color{blue}$\delta$}
	\end{overpic}
	\caption{The shaded region is the decoration $\mathcal{D}$, the red curve is $\beta$, and the blue curve is $\delta$.}\label{fig: decoration_2}
\end{figure}
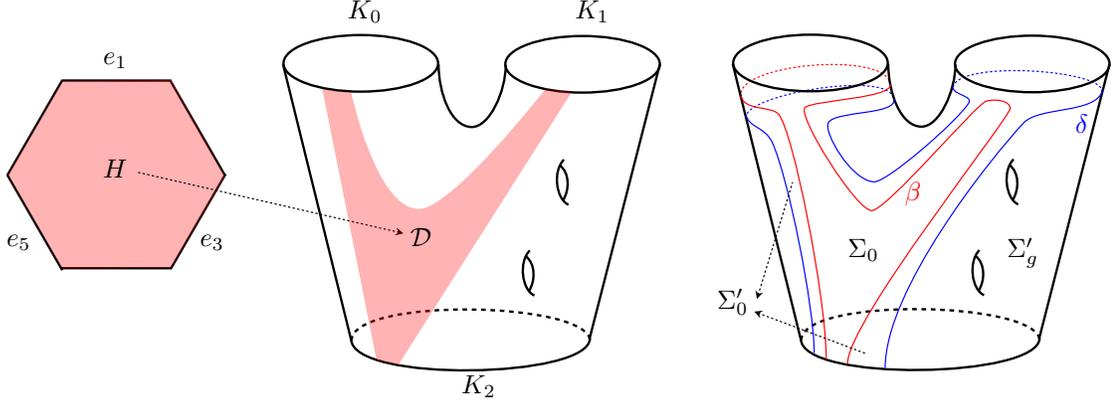

Picking an element $\theta\in KHI^-(Y_1,K_1)$ and taking $\nu_g\in KHI^-(Z_{2g},B_{g})$ as in Definition \ref{defn: the element nu_g}, we can define a cobordism map
$$F_{W,\Sigma_g,\mathcal{D},\theta}:KHI^{-}(Y_0,K_0)\to KHI^{-}(Y_2,K_2)$$
by taking
\begin{equation}\label{eq: F_W,Sigma,D,theta}
	F_{W,\Sigma_g,\mathcal{D},\theta}(x)=F^-_{W-N(\Sigma_g)}\circ\rho'\bigg(\big(\rho(x\otimes\theta)\big)\otimes\nu_g\bigg),
\end{equation}
where
$$\rho:KHI^-(Y_0,K_0)\otimes KHI^-(Y_1,K_1)\to KHI^-(Y_0\# Y_1, K_0\# K_1),$$
and
$$\rho': KHI^-(Y_0\# Y_1, K_0\# K_1)\otimes KHI^{-}(Z_{2g},B_{g})\to KHI^-(Y_0\# Y_1\# Z_{2g}, K_0\# K_1\# B_{g}).$$

Using essentially the same argument as in Section \ref{subsec: two ends}, we can conclude the following.
\blem\label{lem: F_W,Sigma,D,theta is a module morphism}
The map $F_{W,\Sigma_g,\mathcal{D},\theta}$ defined as above is a $\mathbb{C}[U]$-module morphism.
\elem

\begin{lem}\label{lem: functoriality of F_W,Sigma,D,theta}
	Suppose $(W,\Sigma_g,\mathcal{D})$ is a cobordism from $(Y_0,K_0)\sqcup (Y_1,K_1)$ to $(Y_2,K_2)$, where $K_1$ is the unknot. Suppose $(W',\Sigma'_{g'})$ is a cobordism from $\emptyset$ to $(Y_1,K_1)$, as in Section \ref{subsec: one end}. On $\Sigma_{g}\cup\Sigma'_{g'}$, we can smoothly round the corners of $\mathcal{D}$ that are in the interior of $\Sigma_{g}\cup\Sigma'_{g'}$ to make it a decoration in the sense of Section \ref{subsec: two ends}. Then there is an identity
	$$F_{W\cup W',\Sigma_g\cup\Sigma'_{g'},\mathcal{D}}=F_{W,\Sigma_g,\mathcal{D},\theta(W',\Sigma'_{g'})}.$$
\end{lem}

\section{The tube attachment lemma}\label{sec: tube attachment}
\bprop\label{prop: tube attachment lemma}
Let $Y_1=Y_2=S^3$ and $W=[0,1]\times S^3$. Suppose $\Sigma_g\subset W$ is a properly embedded surface so that for $i=0,1$, $\Sigma_g\cap{i}\times S^3$ is a knot $K_i$. Let $\mathcal{D}\subset \Sigma_g$ be a decoration in the sense of Section \ref{subsec: two ends}. Let $\alpha\subset {\rm int}(W)$ be a connected simple arc so that $\alpha\cap\Sigma_g=\partial \alpha$ and $\alpha\cap\mathcal{D}=\emptyset$. We can remove the neighborhood of $\partial \alpha\subset \Sigma_g$ and attach a tube to $\Sigma_g$ along $\alpha$ to create a new surface $\Sigma_{g+1}\subset W$. Then $\mathcal{D}$ is still a decoration on $\Sigma_{g+1}$. Then, we have the following identity. 
$$F_{W,\Sigma_{g+1},\mathcal{D}}=U\cdot F_{W,\Sigma_g,\mathcal{D}}.$$
\eprop

\bpf
Let $B\subset \int(\Sigma_g)$ be a disk that contains $\partial \al$ and $B\cap\mathcal{D}=\partial B\cap \partial \mathcal{D}\subset \eta([0,1]\times\{1\}).$ Let $N_1$ be a $4$-ball neighborhood of $B\subset W$ and $N_2$ be a $4$-ball neighborhood of $\al\subset W$ so that $W''=N_1\cup N_2\cong S^1\times D^3$. Let $W'=W-W''$, ${\Sigma}'_{g}=\Sigma_{g+1}\cap W'$, $\Sigma''_0=\Sigma_g\cap W''$, and $\Sigma''_1=\Sigma_{g+1}\cap W''.$ Let $\mathcal{D}'=\mathcal{D}\cap \Sigma'_{g}$. Note that $\partial W''\cong S^1\times S^2$, 
$$\mathbb{U}_1=\Sigma''_0\cap \partial W''= \Sigma''_1\cap \partial W''$$
is an unknot, and $\mathcal{D}'$ is a decoration on $\Sigma'_g$ in the sense of Section \ref{subsec: three ends}. By Lemma \ref{lem: functoriality of F_W,Sigma,D,theta}, we know
$$F_{W,\Sigma_g,\mathcal{D}}=F_{W',\Sigma'_g,\mathcal{D}',\theta(W'',\Sigma''_0)},$$
and
$$F_{W,\Sigma_{g+1},\mathcal{D}}=F_{W',\Sigma'_g,\mathcal{D}',\theta(W'',\Sigma''_1)}.$$
Hence, the proof of Proposition \ref{prop: tube attachment lemma} is reduced to proposition \ref{prop: U multiplication} below by the fact that the cobordism maps are $\mathbb{C}[U]$-module homomorphisms as in Lemma \ref{lem: F_W,Sigma,D is a module morphism} and Lemma \ref{lem: F_W,Sigma,D,theta is a module morphism}.
\epf
\blem \label{lem : special cobordism for model}
The special cobordism $D^- {N({\Sigma_1''})}$ coincides with the cobordism $X^1_0$ in Section \ref{subsec: borromean knots}.
\elem

\begin{figure}[ht]
	\begin{overpic}[width=0.8\textwidth]{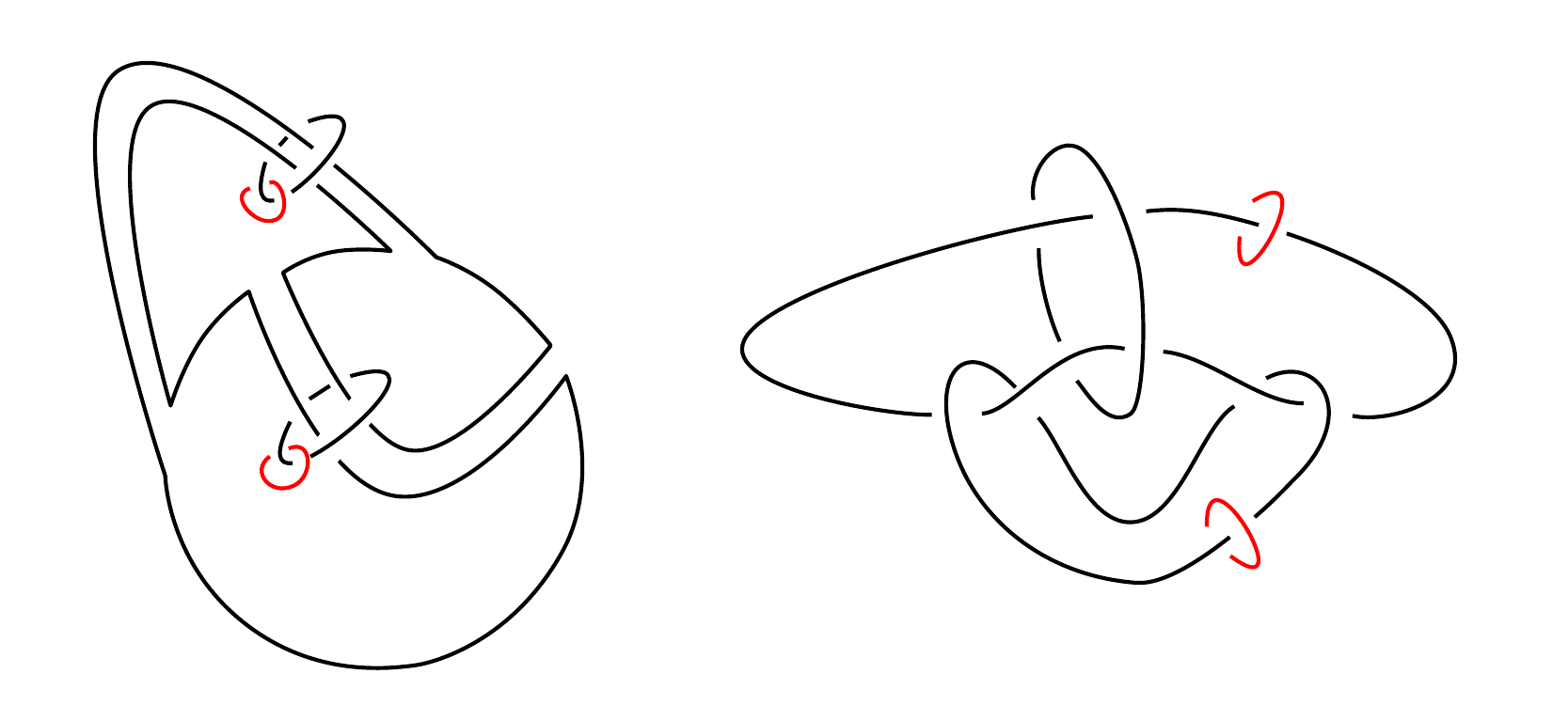}
	\put(19, 38) {0}
    \put(15, 29) {0} 
    \put(23, 22) {0}
    \put(16, 12) {0}
    \put(59, 30) {0}
    \put(59, 21) {0}
    \put(59, 30) {0}
    \put(78, 26) {0}
    \put(78, 13) {0}
	\end{overpic}
	\vspace{0.1in}
	\caption{Left: the special cobordism $D^- {N({\Sigma_1''})}$. Right: the manifold $X_{01}$. The red circles denote the attaching cycles of the 4-dimensional 2-handles for the special cobordism.}\label{fig: Borromean_2}
\end{figure}
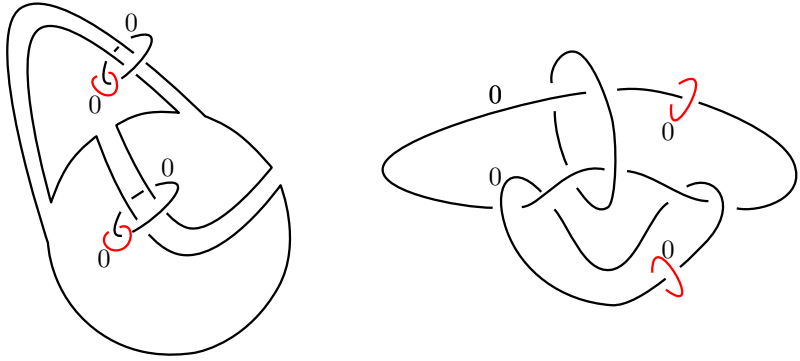

\bpf
We observe that $\Sigma_1''$ can be obtained by attaching two 2-dimensional 1-handles to the unknot which are dual to each other as in Figure \ref{fig: Borromean_2}. We now introduce a cancelling pair of 4-dimensional 1- and 2-handles that define the product $S^3 \times I$ so that each 1-handle is a tubular neighbourhood of the 2-dimensional 1-handle or the band attached to the unknot.  The 1-handles now are contained in $N(\Sigma_1'')$, and therefore disappears when we consider $D^- {N({\Sigma_1''})}$. Hence the special cobordisms are given by two 4-dimensional 2-handles. See Figure \ref{fig: Borromean_2}.
\epf
\bprop\label{prop: U multiplication}
Under the notations in the proof of Proposition \ref{prop: tube attachment lemma}, we have the following identity:
$$\theta(W'',\Sigma''_1)=U\cdot \theta(W'',\Sigma''_0).$$
\eprop

\bpf

In \cite{LY2020}, the authors constructed an $\intg_2$-grading on $KHI^-$, such that the $U$ action preserves this $\intg_2$-grading. We can take $\nu_g$ to be a homogeneous element with respect to this $\intg_2$-grading. Since all the maps involved in the paper come from cobordism maps, they are all homogeneous. Note that $W''$ is diffeomorphic to $S^1\times D^3$, so we can attach a $4$-dimensional $2$-handle to make it a $D^4$. Also, $\partial W''\cong S^1\times S^2$, and this $4$-dimensional $2$-handle corresponds to a Dehn surgery along $S^1$ to make it $S^3$. Since this $4$-dimensional handle attachment is disjoint from the knot, a routine argument shows that it induces a grading-preserving map between $KHI^-$, which commutes with maps coming from contact handle attachments along the boundary of the knot complement. As a result, we have the following commutative diagram:
\begin{equation*}
	\xymatrix{
	&&KHI^-(Z_2,B_1)\ar[rr]^{=}\ar[d]^{F^-_{W''-N(\Sigma_1'')}}&&KHI^-(Z_2,B_1)\ar[d]^{F^-_{D^4-N(\Sigma_1'')}}\\
	&&KHI^-(Z_1,\mathbb{U}_1)\ar[rr]^{F_1^-}\ar[d]^{F^-_{\mathbb{U}_1}}&&KHI^-(S^3,\mathbb{U}_1)\ar[d]^{F^-_{\mathbb{U}_1}}\\
	I^{\#}(S^3)\ar[rr]^{G^{\#}_1}&&I^{\#}(Z_1)\ar[rr]^{F^{\#}_1}&&I^{\#}(S^3)
	}
\end{equation*}
Here $F_1^-$ and $F_1^{\#}$ are the maps associated to the $4$-dimensional $2$-handle attachment, and $G^{\#}_1$ is obtained by attaching a $4$-dimensional $1$-handle to $B^4$ that forms $W''$. The $4$-dimensional $1$- and $2$-handles forms a canceling pair, so we know
$$F_1^{\#}\circ G_1^{\#}=id.$$
As a result $F_1^{\#}$ is surjective. On the other hand, as above it correspond to a surgery along $S^1\subset S^1\times S^2$, so it kills the class $[S^1]\in H_1(S^1\times S^2)$. As a result, using $H_1$-actions, we can show that
$${\rm ker}(F_1^{\#})={\rm ker}\mu([S^1])={\rm Span}(1)\subset I^{\#}(S^1\times S^2).$$
(Recall $I^{\#}(S^1\times S^2)$ has two generators $1$ and $\theta_1$ as in Section \ref{subsec: borromean knots}) Hence we must have $F_1^{\#}(\theta_1)=1\in I^{\#}(S^3)$. 

Note that we start with the element $\nu_1\in KHI^-(Z_2,B_1)$. From Lemma \ref{lem : special cobordism for model}, we know that the special cobordism $D^-N(\Sigma_1'')$ coincides with the cobordism $X^1_0$ in Section \ref{subsec: borromean knots}, so we conclude from Lemma \ref{lem: the map F^2n_0} that 
$$F^-_{D^4-N(\Sigma_1'')}(\nu_1)=U\in KHI^-(S^3,\mathbb{U}_1,-1).$$
Furthermore, we know that
$$F^-_{\mathbb{U}_1}\circ F^-_{D^4-N(\Sigma_1'')}(\nu_1)=1\in I^{\sharp}(S^3).$$
Since $F^-_1$ is grading-preserving as discussed above, we know
$$F^-_{W''-N(\Sigma''_1)}(\nu_1)\in KHI^-(S^1\times S^2,\mathbb{U}_1,-1).$$
Note that by Lemma \ref{lem: minus version for the unknot} we have
$$KHI^-(S^1\times S^2,\mathbb{U}_1,-1)={\rm Span}(U\cdot 1,U\cdot \theta_1)$$
and we have an isomorphism
$$F^-_{\mathbb{U}_1}: KHI^-(S^1\times S^2,\mathbb{U}_1,-1)\to I^{\#}(S^1\times S^2)$$
Since all elements are homogeneous, the only possibility is that
$$\theta(W'',\Sigma''_1)=F^-_{W''-N(\Sigma''_1)}(\nu_1)=U\cdot \theta_1.$$

The computation for $\theta(W'',\Sigma''_0)$ is essentially the same, but simpler as after attaching the $4$-dimensional $2$-handle, the special cobordism $D^4-N(\Sigma''_0)$ is simply a product one. As a result, we can conclude that
$$\theta(W'',\Sigma''_0)=\theta_1\in KHI^{-}(S^1\times S^2, \mathbb{U},0).$$
and hence
$$\theta(W'',\Sigma''_1)=U\cdot\theta(W'',\Sigma''_0).$$
\epf

\section{Torsion order and connected sum}
In this section, we prove the following result.
\bprop\label{prop: torsion order and connected sum}
Suppose $K_1$ and $K_2$ are two knots in $S^3$. Then we have the following equality
$${\rm ord}_{U}(K_1\# K_2)=\max\{{\rm ord}_{U}(K_1),{\rm ord}_{U}(K_2)\}.$$
\eprop

\blem
Suppose $A$ and $B$ are two $\intg$-graded finitely generated rank-one $\mathbb{C}[U]$-modules so that $U$ has degree $-1$. Then
$${\rm ord}_{U}(A\widetilde{\otimes}B)=\max\{{\rm ord}_{U}(A),{\rm ord}_{U}(B)\}.$$
\elem
\bpf

We denote the $U$-action on $A$ as $U_1$, and the action on $B$ as $U_2$. Choose $a\in A$ as a $U_1$-free element at the highest grading in $A$, and $b\in B$ as a $U_2$-free element at the highest grading in $B$. We can treat $A$ and $B$ as complex vector spaces and select a set of homogeneous bases for each, ensuring that the following conditions are met.
\begin{enumerate}
	\item The elements $U_1^{j}a$ for $j\in\intg_{\geq 0}$ are all in the basis for $A$, and any other element in the basis of $A$ is $U$-torsion.
	\item The elements $U_2^{j}b$ for $j\in\intg_{\geq 0}$ are all in the basis for $B$, and any other element in the basis of $B$ is $U$-torsion.
\end{enumerate}
Assume, without loss of generality, we have
$$r={\rm ord}_{U}(A)\geq {\rm ord}_{U}(B).$$

{\bf Step 1.} We show that ${\rm ord}_{U}(A\widetilde{\otimes}B)\geq \max\{{\rm ord}_{U}(A),{\rm ord}_{U}(B)\}.$ 
Take $a'\in A$ a homogeneous element so that $U^r_1a'\neq0$ but $U^{r+1}_1a'=0$. Now take the element $x=a'\otimes b\in A\otimes B$. We would like to show that $U_1^rx\notin \im(U_1+U_2)$.

Assume the contrary, {\it i.e.}, we have
\begin{equation*}
	(U_1+U_2)\sum_{i=1}^n\lambda_i\cdot(a_i\otimes b_i)=a'\otimes b.
\end{equation*}
where all $a_i$ and $b_i$ are included in the chosen bases for $A$ and $B$, respectively. Since $b$ is a free element in the basis, there must be at least one free element among the $b_i$. Without loss of generality, we can assume that $b_1$ is the free element with the lowest grading, that $b_2$ through $b_m$ are equal to $b_1$, and that $b_i$ differs from $b_1$ for $i$ ranging from $m+1$ to $n$. It is important to note that $B$ has a rank of one as a $\mathbb{C}[U]$-module, which means $b_1$ is the sole free element in the basis at its particular grading. Consequently, we observe the following:
\begin{equation}\label{eq: a'b in image, 2}
	a'\otimes b=\sum_{i=1}^m\lambda_i\cdot a_i\otimes (U_2b_1)+{\rm (other~higher~grading~or~}U_2\dash{\rm torsion~terms)}.
\end{equation}
Thus, we conclude that $m=1$, $\lambda_1=1$, $a'=a_1$, and $b=U_2b_1$. However, this conclusion contradicts the selection of $b$ as a $U$-free element of maximal grading.

Consequently, $U^r_1x$ does not fall within the image of $U_1-U_2$, indicating that $[U_1^rx]$ is non-zero in ${\rm coker}(U_1-U_2)$ and subsequently embeds into $A\widetilde{\otimes}B$. Nonetheless, $U_1^{r+1}x = (U^{r+1}_1a' \otimes b) = 0$. Therefore, we deduce that $[x]$ corresponds to a $U$-torsion element in $A\widetilde{\otimes}B$ with order $r$, leading us to the following conclusion
$${\rm ord}_{U}(A\widetilde{\otimes}B)\geq \max\{{\rm ord}_{U}(A),{\rm ord}_{U}(B)\}.$$

{\bf Step 2.} We show that ${\rm ord}_{U}(A\widetilde{\otimes}B)\leq \max\{{\rm ord}_{U}(A),{\rm ord}_{U}(B)\}.$ 
We can regard $A\widetilde{\otimes}B$ as $\ker (U_1-U_2)\oplus{\rm coker}(U_1-U_2)$. As a result, we treat the kernel and co-kernal parts separately. First, suppose
$$x=\sum_{i=1}^n\lambda_i\cdot(a_i\otimes b_i)\in\ker (U_1-U_2).$$
If we take $a'\otimes b=0$ in Equation (\ref{eq: a'b in image, 2}), we conclude that no $b_i$ could be $U_2$-free. As a result, we conclude that any such element has torsion order at most ${\rm }ord_{U}(B)$.

Second, suppose
$$x=\sum_{i=1}^n\lambda_i\cdot (a_i\otimes b_i)\notin \im(U_1-U_2)$$
is an element so that $[x]\in {\rm coker}(U_1+U_2)$ is $U$-torsion. We may assume that for $i=1,...,m$, both $a_i$ and $b_i$ are free, and for $i=m+1,...,n$, either $a_i$ or $b_i$ is torsion. Take $k$ large enough, we conclude that
$$U_1^kU_2^kx=\sum_{i=1}^m\lambda_i\cdot (U_1^ka_i)\otimes(U_2^kb_i).$$
Note that at sufficiently low gradings, the torsion components of both $A$ and $B$ do not contribute to $A\widetilde{\otimes}B$. Therefore, the low-grade portions of $A\widetilde{\otimes}B$ correspond to the standard model of $\mathbb{C}[U]\widetilde{\otimes}\mathbb{C}[U]$. Consequently, it is evident that if all $a_i$ and $b_i$ are free, and $k$ is sufficiently large, then:
$$\sum_{i=1}^m\lambda_i\cdot (U_1^ka_i)\otimes(U_2^kb_i)\in\im(U_1+U_2)\Leftrightarrow \sum_{i=1}^m\lambda_i\cdot (a_i)\otimes(b_i)\in\im(U_1+U_2).$$
As a result, we conclude that
$$[x]=[y=\sum_{i=m+1}^n\lambda_i\cdot (a_i\otimes b_i)]\in {\rm coker}(U_1+U_2).$$
Now for each $i=m+1,...,n$, either $a_i$ or $b_i$ is torsion, so we conclude that $[x]$ has torsion order at most $r={\rm ord}_{U}(A)\geq {\rm ord}_{U}(B)$.
\epf

\section{Bridge index, torsion order, and Dehn surgeries}
\bpf[Proof of Theorem \ref{thm: torsion order and bridge index}]
Suppose $K\subset S^3$ is a knot with bridge index $b$. A well-known observation is the following ({\it c.f.} \cite[Section 1.6]{juhasz2020bridge}). 
\blem\label{lem: b band to unlink}
For a knot $K$ of bridge index $b$, we can attach $(b-1)$ bands to $K\#\widebar{K}$ to make it a $b$-component unlink $\mathbb{U}_{b}$.
\elem
Consider the cobordism $S$ arising from Lemma \ref{lem: b band to unlink} in the manifold $W=[0,1]\times S^3$, connecting $K\#\widebar{K}\subset {0}\times S^3$ to $\mathbb{U}_{b}\subset {1}\times S^3$. Let $\widebar{S}$ be the horizontal mirror image of $S$. As demonstrated in the proof of \cite[Lemma 4.1]{juhasz2020bridge}, the cobordism $S\cup \widebar{S}$ can be derived from a concordance by attaching $(b-1)$ tubes in two different ways:
\begin{enumerate}
	\item $S\cup\widebar{S}$ can be obtained from a concordance $R\cup \widebar{R}$ by attaching $(b-1)$ tubes, where $R$ is obtained from $S$ by capping off $(b-1)$ components of the unlink $\mathbb{U}_b$.
	\item $S\cup\widebar{S}$ can be obtained from a product concordance $P$ by attaching $(b-1)$ tubes.
\end{enumerate}

Essentially, when considering $R\cup\widebar{R}$, attaching tubes is the reverse process of capping off. To understand why there is a product cobordism, observe that $S$ is derived from the product concordance by attaching $(b-1)$ bands. Each band, along with its mirror image in $\widebar{S}$, results in a tube being attached to the product cobordism. Let $\mathbf{1}\in KHI^{-}(Y_{2g},B_g,g)$ be the generator for any $g$. We can choose a decoration $\mathcal{D}$ that is independent of any of these constructions. Now, by applying Proposition \ref{prop: tube attachment lemma}, we conclude that:
$$U^{b-1}\cdot F_{W, R\cup \widebar{R}, \mathcal{D}}=F_{W, S\cup \bar{S},\mathcal{D}}=U^{b-1}\cdot F_{W,P,\mathcal{D}}:KHI^{-}(K\#\widebar{K})\ra KHI^{-}(K\#\widebar{K}).$$
Note that $R\cup \widebar{R}$ factors through the unknot $\mathbb{U}_1$. As a result, we can conclude, using Lemma \ref{lem: minus version for the unknot} and Lemma \ref{lem: rho is an iso}, that $U^{b-1}=0$ on the torsion part of $KHI^{-}(K\#\widebar{K})$. This implies that:
$${\rm ord}_U(K\#\widebar{K})\leq b-1={\rm br}(K)-1.$$
This completes the proof of Theorem \ref{thm: torsion order and bridge index} as per Proposition \ref{prop: torsion order and connected sum}.
\epf

\bpf[Proof of Theorem \ref{thm: main}]
Suppose $K\subset S^3$ is an alternating knot of bridge index $\leq 3$, then Theorem \ref{thm: torsion order and bridge index} and Lemma \ref{lem: odd torsion order} implies that ${\rm ord}_U(K)=1$ and hence Proposition \ref{prop: dehn surgery on torsion order one} applies and we are done.
\epf

%\tableofcontents%table of contents
% \newpage

%————Start from here————

%————End from here————
%\newpage
\bibliographystyle{alpha}
\bibliography{ref.bib}

\end{document}